\documentclass[a4paper,11pt]{article}
\usepackage[centertags]{amsmath}
\usepackage{amssymb}
\usepackage{amsthm}
\usepackage{graphicx}
\usepackage{amsfonts}
\usepackage{mathrsfs}
\usepackage{times,cite}
\usepackage{geometry}
\usepackage{color}
\newtheorem{thm}{Theorem}[section]
\newtheorem{lem}[thm]{Lemma}
\newtheorem{prop}[thm]{Proposition}
\newtheorem{cor}[thm]{Corollary}

\theoremstyle{definition}

\newcommand{\e}{\varepsilon}
\newcommand{\p}{\partial}
\newcommand{\y}{\langle y \rangle}
\newcommand{\x}{\langle \eta \rangle}
\newenvironment {Proof} {\noindent {\bf Proof. }}{\hfill $\square$\par\vspace{3mm}}
\newenvironment {ProofTheorem} {\noindent {\bf Proof of Theorem 1.1.}}{\hfill $\square$\par\vspace{3mm}}
\numberwithin{equation}{section}
\allowdisplaybreaks    

\newenvironment{acknowledgment}{\noindent{\bf Acknowledgment }}{}
\date{}
\title{Symmetrical Prandtl boundary layer expansions of steady Navier-Stokes equations on bounded domain}
\author{{ Shijin Ding$^a$, Quanrong Li$^b$\thanks{Corresponding author. Emails: {\it dingsj@scnu.edu.cn(S. Ding), quanrong\_li@szu.edu.cn(Q. Li)}}}\\
{\it\small $^a$South China Research Center for Applied Mathematics and Interdisciplinary Studies}\\ {\it\small South China Normal University}\\
{\it\small Guangzhou, 510631, Guangdong, China}\\
{\it\small $^b$College of Mathematics and Statistics}\\ {\it\small Shenzhen University}\\
{\it\small Shenzhen, 518060, Guangdong, China}\\}
\begin{document}
\newcommand{\D}{\displaystyle}
\maketitle

 \begin{abstract}
 This paper is concerned with the validity of the Prandtl boundary layer theory in the inviscid limit of the steady incompressible Navier-Stokes equations, which is an extension of the pioneer paper \cite{GN14}(Y. Guo et al., 2017, Ann. PDE) from a domain of $[0,L]\times\mathbb{R}_+$ to $[0,L]\times[0,2]$. Under the symmetry assumption, we establish the validity of the Prandtl boundary layer expansions and the error estimates. The convergence rate as $\e\rightarrow 0$ is also given.
 \end{abstract}


\section{Introduction}
\subsection{Formulation of the problem}

In this paper, we consider the following steady incompressible Navier-Stokes
equations
\begin{align}
\begin{cases}\label{1.1}
UU_X+VU_Y+P_X=\e U_{XX}+\e U_{YY},\\
UV_X+VV_Y+P_Y=\e V_{XX}+\e V_{YY},\\
U_X+V_Y=0,
\end{cases}
\end{align}
in the domain
\[\Omega:=\{(X,Y)|0\leq X\leq L, 0\leq Y\leq 2\}\]
 with moving boundary conditions
\[U(X,0)=U(X,0)=u_b>0,~V(X,0)=V(X,2)=0.\]

We will focus on the problem when $\e\rightarrow 0$. As $\e\rightarrow 0$, a formal limit of the solution of (\ref{1.1}) should be the shear flow $[U_0,V_0]=[u_e^0(Y),0]$, which satisfies the corresponding Euler equations. We assume that this smooth positive function $u^0_e(\cdot)$ satisfies $u_e^0(1-Y)=u_e^0(1+Y),$ for any $Y\in[0,1]$ and $u_e^0(0)=u_e^0(2)=u_e\neq u_b$. Accordingly, we assume that the solution $[U,V]$ to (\ref{1.1}) satisfies the following symmetrical conditions with respect to $Y=1$
\begin{align*}
U(X,1-Y)=U(X,1+Y),V(X,1-Y)=-V(X,1+Y), Y\in[0,1].
\end{align*}

It should be noted that, due to this assumption, the pair $[U,V]_{1\leq Y\leq 2}$ satisfies equations (\ref{1.1}) as long as  $[U,V]_{0\leq Y\leq 1}$ does. Then our discussion can be restricted to the domain
\[\Omega_0:=\{(X,Y)|0\leq X\leq L, 0\leq Y\leq 1\}\]
and the boundary conditions turn to
\[[U,V](X,0)=[u_b,0],\quad [U_Y,V](X,1)=[0,0].\]

Now we introduce the Prandtl's scaling
\[x=X,y=\frac{Y}{\sqrt\e}\]
and the new unknown functions
\[U^\e(x,y)=U(X,Y),V^\e(x,y)=\frac{1}{\sqrt{\e}}V(X,Y).\]
Under this transformation, system (\ref{1.1}) can be rewritten as
\begin{align}\label{1.2}
\begin{cases}
U^\e U_x^\e+V^\e U_y^\e+P_x^\e=U_{yy}^\e+\e U_{xx}^\e,\\
U^\e V_x^\e+V^\e V_y^\e+P_y^\e/\e=V_{yy}^\e+\e V_{xx}^\e,\\
U_x^\e+V_y^\e=0,
\end{cases}
\end{align}
in the domain
\[\Omega_\e:=\{(x,y)|0\leq x\leq L,0\leq y\leq \frac{1}{\sqrt\e}\},\]
with the boundary conditions
\begin{align}\label{1.2a}
[U^\e,V^\e](x,0)=[u_b,0],\quad [U_y^\e,V^\e](x,\frac{1}{\sqrt\e})=[0,0].
\end{align}

In what follows, we intend to find the exact solutions $[U^\e,V^\e,P^\e]$ in form of
\begin{align}\label{1.3}
\begin{cases}
U^\e(x,y)=u_{app}(x,y)+\e^{\gamma+\frac{1}{2}}u^\e(x,y),\\
V^\e(x,y)=v_{app}(x,y)+\e^{\gamma+\frac{1}{2}}v^\e(x,y),\\
P^\e(x,y)=p_{app}(x,y)+\e^{\gamma+\frac{1}{2}}p^\e(x,y),
\end{cases}
\end{align}
where
\begin{align}\label{1.4}
\begin{cases}
u_{app}(x,y)=u_e^0(\sqrt\e y)+u_p^0(x,y)+\sqrt\e u_e^1(x,\sqrt\e y)+\sqrt\e u^1_p(x,y),\\
v_{app}(x,y)=v_p^0(x,y)+v_e^1(x,\sqrt\e y)+\sqrt\e v^1_p(x,y),\\
p_{app}(x,y)=\sqrt\e p_e^1(x,\sqrt\e y)+\sqrt\e p^1_p(x,y)+\e p_p^2(x,y).
\end{cases}
\end{align}

Substituting (\ref{1.3}) into (\ref{1.2}), we get
\begin{align}\label{1.5}
\begin{cases}
R^u_{app}+\e^{\gamma+\frac{1}{2}}\left[(u^\e\p_x+v^\e\p_y)u_{app}+(u_{app}\p_x+v_{app}\p_y)u^\e
+p_x^\e-\Delta_\e u^\e\right]\\
\qquad+\e^{2\gamma+1}(u^\e\p_x+v^\e\p_y)u^\e=0,\\
R^v_{app}+\e^{\gamma+\frac{1}{2}}\left[(u^\e\p_x+v^\e\p_y)v_{app}+(u_{app}\p_x+v_{app}\p_y)v^\e
+p_y^\e/\e-\Delta_\e v^\e\right]\\
\qquad+\e^{2\gamma+1}(u^\e\p_x+v^\e\p_y)v^\e=0,\\
\p_xu_{app}+\p_yv_{app}+\e^{\gamma+\frac{1}{2}}(u^\e_x+v^\e_y)=0,
\end{cases}
\end{align}
where $\Delta_\e:=\p_y^2+\e\p_x^2$ and the errors caused by the approximation
\begin{align*}
&R_{app}^u:=(u_{app}\p_x+v_{app}\p_y)u_{app}+\p_xp_{app}-\Delta_\e u_{app},\\
&R_{app}^v:=(u_{app}\p_x+v_{app}\p_y)v_{app}+\p_yp_{app}/\e-\Delta_\e v_{app},
\end{align*}
or precisely
\begin{align}
R_{app}^u:=&\left[(u_e^0+u_p^0+\sqrt\e[u_e^1+u_p^1])\p_x+(v_p^0+v_e^1+\sqrt\e v_p^1)\p_y\right](u_e^0+u_p^0+\sqrt\e[u_e^1+u_p^1])\nonumber\\
&+\p_x(\sqrt\e[p_e^1+p_p^1]+\e p_p^2)-(\p_y^2+\e\p_x^2)(u_e^0+u_p^0+\sqrt\e[u_e^1+u_p^1]),\label{1.6}\\
R_{app}^v:=&\left[(u_e^0+u_p^0+\sqrt\e[u_e^1+u_p^1])\p_x+(v_p^0+v_e^1+\sqrt\e v_p^1)\p_y\right](v_p^0+v_e^1+\sqrt\e v_p^1)\nonumber\\
&+\p_y(p_e^1+p_p^1+\sqrt\e p_p^2)/\sqrt\e-(\p_y^2+\e\p_x^2) (v_p^0+v_e^1+\sqrt\e v_p^1).\label{1.7}
\end{align}

Now the boundary conditions can be rewritten as
\begin{align}\label{1.8}
\begin{cases}
u_{app}(x,0)+\e^{\gamma+\frac{1}{2}}u^\e(x,0)=u_b,\quad \p_yu_{app}(x,\frac{1}{\sqrt\e})+\e^{\gamma+\frac{1}{2}}u_y^\e(x,\frac{1}{\sqrt\e})=0,\\
v_{app}(x,0)+\e^{\gamma+\frac{1}{2}}v^\e(x,0)=0,\quad v_{app}(x,\frac{1}{\sqrt\e})+\e^{\gamma+\frac{1}{2}}v^\e(x,\frac{1}{\sqrt\e})=0.
\end{cases}
\end{align}

It is clear that there are only three equations with two boundary conditions, but there are twelve unknown functions, which makes this system unclosed. To construct the approximate solution, we have to divide this big system into a few subsystems in terms of the order of $\e$.

\subsection{Boundary conditions}
Let us see how to impose boundary conditions for each subsystem. For convenience, denote $z:=\sqrt\e y$.

Boundary conditions on $\{y=0\}$:
\begin{align}
&u^0_e(0)+u^0_p(x,0)=u_b,\quad~&&u^1_e(x,0)+u^1_p(x,0)=0,\quad~&u^\e(x,0)=0;\label{1.9}\\
&v^0_p(x,0)+v^1_e(x,0)=0,\quad~&&v^1_p(x,0)=0,\quad~&v^\e(x,0)=0.\label{1.10}
\end{align}

Boundary conditions on $\{y=\frac{1}{\sqrt\e}\}$:
\begin{align}
&u^0_{py}(x,\frac{1}{\sqrt\e})=0,\quad~ &u^1_{ez}(x,1)=0,\quad~ &u^1_{py}(x,\frac{1}{\sqrt\e})=0,\quad~ &u^\e_y(x,\frac{1}{\sqrt\e})=0;\label{1.11}\\
&v^0_p(x,\frac{1}{\sqrt\e})=0,\quad~&v^1_e(x,1)=0,\quad~&v^1_p(x,\frac{1}{\sqrt\e})=0,\quad~&v^\e(x,\frac{1}{\sqrt\e})=0.\label{1.12}
\end{align}

Boundary conditions on $\{x=0\}$ :
\begin{align}
&u^0_p(0,y)=\bar{u}_0(y),\quad~&&u^1_e(0,z)=u^1_b(z),\quad~&&u^1_p(0,y)=\bar{u}_1(y),\quad~u^\e(0,y)=0;\label{1.13}\\
&v^1_e(0,z)=V_{b0}(z),\quad~&& v^\e(0,y)=0.&& \label{1.14}
\end{align}

Boundary conditions on $\{x=L\}$:
\begin{align}\label{1.15}
v^1_e(L,z)=V_{bL}(z),\quad~[p^\e-2\e u^\e_x,u_y^\e+\e v^\e_x](L,y)=0.
\end{align}

Denote $u_e:=u^0_e(0)$, which, in general, is not equal to $u_b$. Then by the first condition in (\ref{1.9}), we shall take $u^0_p(x,0)=u_b-u_e$. Similarly, we will take $u^1_p(x,0)=-u^1_e(x,0)$ and $v^1_e(x,0)=-v^0_p(x,0)$, as $u^1_e(x,0)$ and $v^0_p(x,0)$ will be defined automatically by the profile $u^1_e$ and $v^0_p$, respectively.

For the existence of the Euler corrector $[u^1_e,v^1_e,p^1_e]$, it is necessary for us to impose the following compatibility conditions:
\[V_{b0}(0)=-v^0_p(0,0),\quad~V_{bL}(0)=-v^0_p(L,0),\quad~V_{b0}(1)=V_{bL}(1)=0.\]

In addition, as will be seen in Section 3 that $v_{ezz}^1(x,1)=0$ follows directly from the boundary condition $v^1_e(x,1)=0$ and the elliptic equation $v^1_e$ satisfies, we should also set that $V^{\prime\prime}_{b0}(1)=V^{\prime\prime}_{bL}(1)=0.$ Moreover, the boundary condition $u^1_{ez}(x,1)=0$ follows as soon as the compatibility condition $u^1_{bz}(1)=0$ is given, since that
\[u^1_e(x,z)=u^1_b(z)-\int_0^xv^1_{ez}(s,z)ds,\]
which is a natural solution by the divergence-free condition $u^1_{ex}+v^1_{ez}=0.$

Collecting the functions prescribed in (\ref{1.13}) and (\ref{1.14}), precisely, $\bar{u}_0(y), u^1_b(z), \bar{u}_1(y)$ and $V_{b0}(z)$, one yields the following boundary conditions on $\{x=0\}$ for $(U^\e,V^\e)$, which represent the in-flow conditions:
\begin{align}
&U^\e(0,y)=u^0_e(z)+\bar{u}_0(y)+\sqrt\e u^1_b(z)+\sqrt\e \bar{u}_1(y);\label{1.16}\\
&V^\e(0,y)=v^0_p(0,y)+V_{b0}(z)+\sqrt\e v^1_p(0,y).\label{1.17}
\end{align}
Here, we infer that $v^0_p(0,y)$ and $\sqrt\e v^1_p(0,y)$ are unnecessary to be prescribed since they can be determined respectively by the parabolic equations they satisfy.

Finally, the prescribed conditions in (\ref{1.15}) give the out-flow conditions for $(U^\e,V^\e)$, in which only $v^1_e,u^\e$ and $v^\e$ are prescribed as these profiles satisfy elliptic equations. Physically, the out-flow condition for $(u^\e,v^\e)$ in (\ref{1.15}) is called the stress-free condition.

\subsection{Main result and discussion}
We state our main result of the present paper as follow:
\begin{thm}\label{thm1}
Let $u_b>0$ be a constant tangential velocity of the Navier-Stokes flow on the boundary $\{Y=0\}$, and let $u^0_e(Y)$ be a smooth positive Euler flow satisfies $u^0_{ez}(1)=0$. Suppose that the boundary conditions prescribed in (\ref{1.9})-(\ref{1.15}) hold and compatibility conditions discussed after those boundary conditions in subsection 1.2 are valid. Suppose further that the positive condition $\min_y\{u^0_e(\sqrt\e y)+\bar{u}_0(y)\}>0$ holds. Then there exists a constant $L_0>0$, which depends only on the prescribed data, such that for $0<L\leq L_0$ and $\gamma\in(0,\frac{1}{5})$, the asymptotic expansions stated in (\ref{1.3})-(\ref{1.4}) is a solution to equations (\ref{1.2}) on $\Omega_\e$ together with the corresponding boundary conditions. The approximate solutions appearing in the expansions are constructed in Section 2, 3, 4 and 5, in which the remainder solutions $[u^\e,v^\e]$ satisfies the estimate
\begin{align}\label{1.18}
\|\nabla_\e u^\e\|_{L^2(\Omega_\e)}+\|\nabla_\e v^\e\|_{L^2(\Omega_\e)}+\|u^\e\|_{L^\infty(\Omega_\e)}+\sqrt\e\|v^\e\|_{L^\infty(\Omega_\e)}\leq C_0.
\end{align}
\end{thm}

With this Theorem and the corresponding estimates for each component of the expansions, we can obtain the convergence rate of this sequence as $\e\rightarrow 0$, which indicates the validity of the asymptotic expansions (\ref{1.3})-(\ref{1.4}). Precisely, we have the following
\begin{cor}\label{cor1}
Under the assumptions of Theorem \ref{thm1}, there is an exact solution $[U,V]$ to the original system (\ref{1.1}) on the rectangle domain $[0, L]\times[0, 2]$ with the corresponding boundary conditions, such that
\begin{align}
&\sup_{(X,Y)\in\Omega}\left|U(X,Y)-u^0_e(Y)-u^0_p\left(X,\frac{Y}{\sqrt\e}\right)-\sqrt\e u^1_e(X,Y)\right|\lesssim\e^{\frac{1}{2}};\label{1.19}\\
&\sup_{(X,Y)\in\Omega}\left|V(X,Y)-\sqrt\e v^0_p\left(X,\frac{Y}{\sqrt\e}\right)-\sqrt\e v^1_e(X,Y)\right|\lesssim \e^{\frac{1}{2}+\gamma},\label{1.20}
\end{align}
as $\e\rightarrow 0$, where the zeroth order Prandtl profile $[u^0_p, v^0_p]$ and the first order Euler corrector $[u^1_e,v^1_e]$ are constructed in Section 2 and Section 3, respectively. In particular, in the zero viscosity limit, the convergence $[U,V]\rightarrow [u^0_e,0]$ discussed at the beginning of this paper is valid in the usual $L^p$ norm with convergence rate of order $\e^{\frac{1}{2p}}$, $1\leq p<+\infty$.
\end{cor}

Before continuing, let us give a short historical review on the study of the Prandtl boundary layer theory. It is well known that the Prandtl boundary layer theory was first proposed by L. Prandtl in 1904 in the celebrate lecture `'On fluid motion with very small fraction'' at the Heidelberg mathematical congress, see \cite{Pr1905}. In this lecture, Prandtl used theoretical approach with some simple experiments to show that the flow past a body can be divided into two regions: a very thin layer close to the boundary where the viscosity is important, and the remaining region outside this layer where the viscosity can be neglected. Over more than one hundred years, great achievements have been made on the application of computational fluid mechanics and simulation. However, the rigorous proof for the validity of this theory, at least in general cases, is still uncompleted.

One of the main problem on the road to the validity of the Prandtl boundary layer theory is the well-posedness of the Prandtl equation, which was initiated by O. Oleinik in \cite{Ol63} with $p_x\leq 0$ for the steady setting, and in \cite{Ol67} with assuming monotonic-in-$y$ to the initial data of tangential velocity for the unsteady setting, see also the book \cite{OS99}. Subsequently, these problems attracted considerable attention of many excellent mathematicians. In the steady case, if $p_x>0$, then boundary layer separation will appear in the physical pointview, which has been studied by Goldstein and Stewartson \cite{Gs48,Sk58}, see also \cite{DM15}. For the unsteady case, the local well-posedness of Prandtl equation in $[0,L]\times\mathbb{R}_+$, and global well-posedness for $L$ sufficiently small were obtained in \cite{Ol67,OS99}, by the Crocco transformation. Afterwards, still by the Crocco transformation, this global well-posedness was extended to arbitrary $L<+\infty$ in the sense of weak solution, under the assumption of $p_x\leq 0$ by Z. Xin et al\cite{XZ04}. Without the Crocco transformation, the local well-posedness was also established in \cite{AWXY15,MW15} by energy method under the same monotonicity assumption. So far, the global existence of regular solutions to Prandtl equation is still open, even with the monotonicity assumption. When the monotonicity assumption is generalized to multiple monotonicity regions, the local well-posedness is also valid in the analytic setting\cite{KMVW14}. In the direction of removing the monotonicity assumption, we refer to \cite{SC98a,SC98b,LCS03,KV13,IV16,GVM13} for some results in analytic or Gevery setting, while in the Sobolev setting, the equations are ill-posed(Cf.\cite{GVD10,GVN12}). There are also some results on the finite-time blow-up solutions, see \cite{EE97,KVW15,HH03}.

The main purpose of the present paper is to study the validity of the expansions (\ref{1.3})-(\ref{1.4}) to the solutions of the steady Navier-Stokes equations. In the unsteady cases, the local validity is given by \cite{SC98a,SC98b} in analytic setting, by \cite{GVMM16} with Gevery setting, and by \cite{Mae14} under the assumption  that the initial vorticity distribution is supported away from the boundary, also see \cite{As91,MT08} for other related results. In addition, there are also some proofs for the invalidity in Sobolev spaces, see \cite{Gre00,GGN15a,GGN15b,GGN15c,GN11}. The first study of the validity for the steady incompressible Navier-Stokes equations was due to the pioneer paper by Y. Guo and T. Nguyen \cite{GN14} in which the problem was set on an infinite domain $[0,L]\times\mathbb{R}_+$ with $L$ small, and the limit is a shear Euler flow. Subsequently, S. Iyer extended $L$ to $\infty$ with the constant limit flow $(1,0)$\cite{IY16}. He also obtained the validity result in the case when the limit Euler flow is a non-shear one with $L$ smll \cite{IY17}. Similar result in a rotating disk $[0,\theta_0]\times[R_0,+\infty)$ with $\theta_0$ small is given in \cite{Iy15}.

This paper aims to extend the results of \cite{GN14} to a bounded domain, $y\in [0,2]$, which is more suitable to the physical reality. To our knowledge, so far, there is no results on a rectangle domain. The main difference between this paper and \cite{GN14} is that the boundary layer consist of two components, $\{Y=0\}$ and $\{Y=2\}$, while in \cite{GN14} there is only one component, $\{Y=0\}$. The extra boundary $\{Y=2\}$ makes it difficult to couple with each other in the analysis of the boundary layers. To overcome this difficulty, we assume that the limit Euler flow is symmetrical, i.e. $u^0_e(Y)=u^0_e(2-Y)$, and make effort to construct the symmetrical Prandtl layer expansions. Since the boundary conditions on $\{Y=1\}$ are generated automatically by the symmetry assumptions, we have to deal with them carefully in the construction of each layers.

The detailed novelties of this paper, we think, can be stated in the following comments.

(a) In the step of constructing the zeroth order Prandtl profiles $[u^0_p,v^0_p]$, we first consider to solve the Prandtl equations in $[0,L]\times\mathbb{R}_+$ in order to use the Von Mises transformation. After the solutions are constructed, we construct $[u^0_p,v^0_p]$ in $[0,L]\times I_\e$ by cut-off method which will yield some new error terms and give rise to some new estimates.

(b) The construction of the first order Euler corrector $[u^1_e,v^1_e,p^1_e]$ is done directly on $[0,L]\times[0,1]$, where a trouble boundary term $u^0_e(1)v^1_{ez}(x,1)$  appears. To deal with it, we add the $x$-depending term $\int_0^xu^0_e(1)v^1_{ez}(s,1)ds$ to the pressure $p^1_e$ so that an elliptic equation will be derived for $v^1_e$.

(c) Similar to the idea of constructing $[u^0_p,v^0_p]$, the extension and cutoff to the domain is also used in the construction of the first order Prandtl corrector $[u^1_p, v^1_p, p^1_p]$, where the proof of some weighted estimates is the most difficult part, especially in dealing with $v_{px},v_{pxy}$ and $v_{pxx}$. In \cite{GN14}, the authors stated the result of the weighted estimates $\|\y^nv_{pxy}\|_{L^2_xL^2_y}$ and $\|\y^nv_{pxyy}\|_{L^2_xL^2_y}$ and proved the unweighted one ($n=0$), gave an idea for the proof of the case $n\not=0$ without details which says that one can test equation (\ref{4.15}) by $\y^nv_x,\y^nv_{xx}$ to get the weighted estimates. However, we find that this is not a trivial problem. The main problem is that the low order term $\|v_{px}\|_{L^2_xL^2_y}$ can only be controlled by $\|\y^{1+}v_{pxy}\|_{L^2_xL^2_y}$ but not by $\|v_{pxy}\|$, which leads to the failure of the iteration on the index $n$ as stated in \cite{GN14}. To overcome this difficulty, we use different test functions and weights, say $y^nv_{yy}$. We first establish the weighted estimate $\|y^nv_{yyy}\|_{L^2_xL^2_y}$ and $\|y^nv_{xyyy}\|_{L^2_xL^2_y}$, see(\ref{4.161}), for the solution of the linearized equation (\ref{4.15}). The reason we use the weight $y^n$ but not $\y^n$ is that if one uses the weight $\y^n$, then some extra (bad) boundary terms will appear. Fortunately, after proving the solvability of the original equation by the fixed point theorem, with the weighted estimates for $\|y^nv_{pyyy}\|_{L^2_xL^2_y}$ and $\|y^nv_{pxyyy}\|_{L^2_xL^2_y}$, we can recover the $\y^n$-weighted estimates for $v_{pxy}$ and $v_{pxxy}$ by using the stream function and a new defined function. Of cause, the cutoff from $\mathbb{R}_+$ to $I_\e$ will also produce some extra terms.

(d) The construction of the remainders $[u^\e,v^\e,p^\e]$ is based on the linearized results from \cite{GN14}. We use the contraction mapping theory to prove the existence of the remainders with (\ref{1.18}), compared to
$\|\nabla_\e u^\e\|_{L^2(\Omega_\e)}+\|\nabla_\e v^\e\|_{L^2(\Omega_\e)}+\e^{\frac{\gamma}{2}}\|u^\e\|_{L^\infty(\Omega_\e)}
+\e^{\frac{\gamma}{2}+\frac{1}{2}}\|v^\e\|_{L^\infty(\Omega_\e)}\leq C_0,$
with $0\leq\gamma\leq\frac{1}{4}$ in \cite{GN14}. Therefore, the rate of convergence in (\ref{1.20}) is as fast as $\e^{\frac{1}{2}+\frac{1}{5}}$, whereas in \cite{GN14} the fastest rate is $\e^{\frac{1}{2}+\frac{1}{8}}$.

\textbf{Notations.} Throughout this paper, we shall use the following notations. We shall use $\y=\sqrt{y^2+1}$ and denote $I_\e:=[0,\frac{1}{\sqrt\e}]$.  For convenience, we will use $\|\cdot\|_p~(1\leq p\leq +\infty)$, and $\|\cdot\|_{H^k}~(k\geq 1)$, to denote the usual $L^p$ norm and $W^{k,2}$ norm of functions defining on various domains, such as $\Omega_0$, $\Omega_\e$, and sometimes $\mathbb{R}_+$ and $I_\e$, depending on the context. We also denote $C(\cdot)$ a universal constant, which depends on the given data listed in the parenthesis. Occasionally, we write $C$ or use the notation $\lesssim$ in the estimates for simplification. It should be noted that the uniform estimates are always independent of $L$ and $\e$. The smallness of $L$ depends only on the given data, while $\e$ is always taken to be small sufficiently. Denote that $\chi(\cdot)$ is a smooth cut-off function supported in $[0,1]$ with $\chi(0)=1$, $\chi(1)=0$, $\chi^\prime(0)=\chi^\prime(1)=0$.

In the rest of this paper, we will construct the  zeroth order Prandtl profile $[u_p^0,v_p^0,0]$ in Section 2, construct the first order Euler corrector $[u_e^1,v^1_e,p_e^1]$ and $p_p^2$ in Section 3. After constructing the first order Prandtl corrector $[u_p^1,v_p^1,p_p^1]$ in Section 4, we will, finally, prove the existence of the reminder in Section 5, which completes the proof of the main results of the present paper.

\section{The zeroth order Prandtl profile}

In order to construct the zeroth order Prandtl profile $[u_p^0,v_p^0,0]$, we denote
\[R^u_0:=(u_e^0+u_p^0)\p_x(u_e^0+u_p^0)+(v_p^0+v_e^1)\p_y(u_e^0+u_p^0)-\p_y^2(u_e^0+u_p^0).\]
Since the Euler profile is always evaluated at $(x,z)=(x,\sqrt\e y)$, we note that
\begin{align*}
&\p_xu_e^0=0,\ (v_p^0+v_e^1)\p_yu_e^0=\sqrt\e(v_p^0+v_e^1)u_{ez}^0,\ \p_y^2u_e^0=\e u^0_{ezz},\\
&u_e^0u_{px}^0+v_e^1u_{py}^0=u_eu^0_{px}+v_e^1(x,0)u_{py}^0+\sqrt\e y(u_{ez}^0u^0_{px}+v^1_{ez}u^0_{py})+E^0,
\end{align*}
where $u_e=u^0_e(0)$ and
\begin{align}\label{2.1} E^0=\e\int_0^y\int_y^r\left[u_{ezz}^0(\sqrt\e\tau)u_{px}^0(x,y)+v^1_{ezz}(x,\sqrt\e\tau)u^0_{py}(x,y)\right]d\tau dr.
\end{align}

In view of the divergence-free condition, we let
\begin{align}\label{2.2}
\begin{cases}
(u_e+u_p^0)u^0_{px}+(v^0_p+v_e^1(x,0))u^0_{py}-u^0_{pyy}=0,\\
u^0_{px}+v^0_{py}=0.
\end{cases}
\end{align}
Then, the zeroth order error term $R_0^u$ is reduced to
\begin{align}\label{2.3}
R_0^u=\sqrt\e(v^0_p+v_e^1)u^0_{ez}+\sqrt\e y(u^0_{ez}u^0_{px}+v^1_{ez}u^0_{py})+E^0-\e u^0_{ezz}.
\end{align}

Base on (\ref{1.8}), we give the following boundary conditions
\[u^0_p(x,0)=u_b-u_e,\ u^0_{py}(x,\frac{1}{\sqrt\e})=0,\ [v_p^0+v_e^1](x,0)=0.\]
Since that $u^0_{px}+v^0_{py}=0$, $v^0_p$ can be expressed as
\[v^0_p(x,y)=\int_y^{\frac{1}{\sqrt\e}}u^0_{px}(x,\theta)d\theta,\]
and the coefficient $v^0_p+v_e^1(x,0)$ can be rewritten as
\[v^0_p(x,y)-v_p^0(x,0)=\int_0^yv^0_{py}(x,\theta)d\theta=-\int_0^yu^0_{px}(x,\theta)d\theta.\]
Then the system (\ref{2.2}) is reduced to the following nonlinear parabolic system of $u^0_p$:
\begin{align}\label{2.4}
\begin{cases}
(u_e+u_p^0)u^0_{px}-\int_0^yu^0_{px}d\theta u^0_{py}=u^0_{pyy},\ y\in I_\e,\\
u^0_p(x,0)=u_b-u_e,\ u^0_{py}(x,\frac{1}{\sqrt\e})=0,\ u^0_p(0,y)=\bar{u}_0(y).
\end{cases}
\end{align}
First, we extend the domain $I_\e$ to $\mathbb{R}_+$ with $\lim\limits_{y\rightarrow\infty}u^0_p(x,y)=0$ in place of the boundary condition $u^0_{py}(x,\frac{1}{\sqrt\e})=0$. Since we shall cut-off the domain from $\mathbb{R}_+$ to $I_\e$ after establishing the estimates for the solution, we denote here by $[u_p^\infty,v^\infty_p]$, for distinction.

Now, use the von Mises transformation:
\[\eta:=\int^y_0(u_e+u^\infty_p(x,\theta))d\theta,\qquad\mathbf{w}(x,\eta):=u_e+u^\infty_p(x,y(\eta)),\]
The function $\mathbf{w}$ then solves
\[\mathbf{w}_x=(\mathbf{w}\mathbf{w}_\eta)_\eta,\ \ \textrm{in}\ \ \Omega_\infty:=[0,L]\times\mathbb{R}_+,\]
which is a standard one-dimensional porous medium equation and is solvable over $\Omega_\infty$, at least when $L$ is small\cite{WYLZ}. In addition, by the Maximum Principle of the porous medium equation, we have
\begin{align}\label{2.5}
0<c_0:=\min_y\{u_b,u_e,u_e+\bar{u}_0(y)\}\leq \mathbf{w}\leq \max_y\{u_b,u_e,u_e+\bar{u}_0(y)\}:=\bar{c}_0.
\end{align}
Now, it remains to derive the energy estimates. Since $\mathbf{w}$ does not vanish on the boundary, we introduce $w:=\mathbf{w}-u_e-[u_b-u_e]e^{-\eta}.$ Then $w$ satisfies
\begin{equation}\label{2.6}
\begin{cases}
w_x=[\mathbf{w}w_\eta]_\eta-[u_b-u_e][we^{-\eta}]_\eta-F_\eta,\\
w(x,0)=0,\lim\limits_{\eta\rightarrow\infty}w(x,\eta)=0,
\end{cases}
\end{equation}
where $F(\eta):=[u_b-u_e][u_e+[u_e-u_b]e^{-\eta}]e^{-\eta}.$ Clearly, $\x^nF(\cdot)\in W^{k,p}(\mathbb{R}_+)$, for any $k\geq 0$ and $p\in[1,+\infty]$. In what follows, we will give the regularity estimates for unique smooth solution to system (\ref{2.6}).

First, we introduce the following weighted iterative norm:
\begin{align}\label{2.7}
\mathcal{N}_j(x):=\sum_{k=0}^j\sup_{0\leq s\leq x}\int_{\mathbb{R}_+}\x^n|\p_x^kw|^2
+\sum^j_{k=0}\int_0^x\int_{\mathbb{R}_+}\x^n\mathbf{w}|\p^k_xw_\eta|, j\geq 1.
\end{align}
Multiplying (\ref{2.6})$_1$ by $\x^nw$ and integrating by parts over $\mathbb{R}_+$ leads to
\begin{align}\label{2.8}
\frac{1}{2}\frac{d}{dx}\int\x^n|w|^2+\int\x^n|w_\eta|^2\lesssim\int\x^n[|w||w_\eta|+|w||F_\eta|],
\end{align}
where the positive upper and lower bounds of $\mathbf{w}$ have been used. Applying Cauchy's inequality to the right-hand side of (\ref{2.8}) gives
\begin{align}\label{2.9}
\frac{d}{dx}\int\x^n|w|^2+\int\x^n|w_\eta|^2\lesssim\int\x^n|w|^2+\int\x^n|F_\eta|^2,
\end{align}
which together with the Gronwall's inequality implies that
\begin{align}\label{2.10}
\sup_{0\leq s\leq x}\int\x^n|w|^2+\int_0^x\int\x^n|w_\eta|^2\leq C(L)(\mathcal{N}_0(0)+1).
\end{align}
This means $\mathcal{N}_0(x)\leq C(\mathcal{N}_0(0)+1)$, for some constant $C>0$ depends only on $L,u_e,u_b,\bar{u}_0.$

Next, applying $\p_x$ to (\ref{2.6})$_1$ yields
\begin{align}\label{2.11}
w_{xx}=[\mathbf{w}w_{x\eta}]_\eta+[w_xw_\eta]_\eta-[u_b-u_e][w_xe^{-\eta}]_\eta.
\end{align}
Similarly, multiplying (\ref{2.11}) by $\x^nw_x$ and integrating by parts over $\mathbb{R}_+$, we get
\begin{align}\label{2.12}
&\frac{1}{2}\frac{d}{dx}\int\x^n|w_x|^2+\int\x^n|w_{x\eta}|^2\nonumber\\
&\lesssim \int\x^n[|w_x||w_{x\eta}|+|w_\eta||w_x||w_{x\eta}|+|w_x|^2|w_\eta|+|w_x|^2],
\end{align}
integrating which over $[0,x]$, together with using the Cauchy's inequality, leads to
\begin{align}\label{2.13}
\sup_{0\leq s\leq x}\int\x^n|w_x|^2+\int_0^x\int\x^n|w_{x\eta}|^2&\lesssim \mathcal{N}_1(0)+\int_0^x\int\x^n[|w_x|^2+|w_\eta|^2|w_x|^2]\nonumber\\
&\lesssim \mathcal{N}_1(0)+\int_0^x(1+\|w_\eta\|_\infty^2)\int\x^n|w_x|^2.
\end{align}
To bound $\|w_\eta\|_\infty$, due to the equation (\ref{2.6}), we have
\begin{align*}
|w_\eta|\leq \int_\eta^\infty|w_{\eta\eta}|d\eta\lesssim&\int(|w_x|+|w_\eta|^2+|w|+|F_\eta|)\nonumber\\
\lesssim&\left(\int\x^n|w_x|^2\right)^{1/2}+\int\x^n|w_\eta|^2+1.
\end{align*}
Furthermore, multiplying (\ref{2.6})$_1$ by $\x^nw$, integrating by parts over $\mathbb{R}_+$ and using Cauchy's inequality, we yield
\begin{align}\label{2.14}
\int\x^n|w_\eta|^2\lesssim \int\x^n(|w|^2+|w_x|^2+|F_\eta|^2)\lesssim \mathcal{N}_1(x)+1,
\end{align}
which implies that
\begin{align}\label{2.15}
\|w_\eta\|_\infty\lesssim \mathcal{N}_1(x)+1.
\end{align}
Now, substituting (\ref{2.15}) into (\ref{2.13}) yields
\begin{align}\label{2.16}
\sup_{0\leq s\leq x}\int\x^n|w_x|^2+\int_0^x\int\x^n|w_{x\eta}|^2\lesssim \mathcal{N}_1(0)+\int_0^x(1+\mathcal{N}_1(s))^2\int\x^n|w_x|^2,
\end{align}
and hence, it follows from (\ref{2.16}) together with (\ref{2.10}) that
\begin{align}\label{2.17}
\mathcal{N}_1(x)\leq C(\mathcal{N}_1(0)+1)+\int_0^x(\mathcal{N}_1(s))^3ds,
\end{align}
which, by Gronwall's inequality, implies that $\mathcal{N}_1(x)\leq C(\mathcal{N}_1(0)+1)$, for $L$ sufficiently small.

In what follows, we shall prove the general estimate for $\mathcal{N}_j(x)$ by mathematical induction. Assume that there holds
\begin{align}\label{2.18}
\mathcal{N}_k(x)\leq C(\mathcal{N}_k(0)+1),
\end{align}
for some $k\geq 1$. Then applying $\p_x^{k+1}$ to (\ref{2.6})$_1$, we get
\begin{align}\label{2.19}
\p_x^{k+1}w_x=[\mathbf{w}\p_x^{k+1}w_\eta]_\eta+\sum_{i=0}^kC_{k+1}^i[\p_x^{k+1-i}w\p_x^iw_\eta]_\eta
-[u_b-u_e][\p_x^{k+1}we^{-\eta}]_\eta.
\end{align}
Similarly as above, multiplying it by $\x^n\p_x^{k+1}w$ and integrating over $\mathbb{R}_+$ leads to
\begin{align}\label{2.20}
&\frac{d}{dx}\int\x^n|\p_x^{k+1}w|^2+\int\x^n|\p_x^{k+1}w_\eta|^2\nonumber\\
&\lesssim \int\x^n|\p_x^{k+1}w_\eta||\p_x^{k+1}w|+\int\x^n\sum_{i=0}^k|\p_x^{k+1-i}w||\p_x^iw_\eta||\p_x^{k+1}w_\eta|\nonumber\\
&\quad+\int\x^n\sum_{i=0}^k|\p_x^{k+1-i}w||\p_x^iw_\eta||\p_x^{k+1}w|+\int\x^n|\p_x^{k+1}w|^2,
\end{align}
where the positive upper and lower bounds of $\mathbf{w}$ have been used.

It follows by Cauchy's inequality that
\begin{align}
&\int\x^n|\p_x^{k+1}w_\eta||\p_x^{k+1}w|\leq \delta \int\x^n|\p_x^{k+1}w_\eta|^2+C\int\x^n|\p_x^{k+1}w|^2,\label{2.21}\\
&\int\x^n\sum_{i=1}^k|\p_x^{k+1-i}w||\p_x^iw_\eta||\p_x^{k+1}w_\eta|\nonumber\\
&\leq \delta \int\x^n|\p_x^{k+1}w_\eta|^2+C\sum_{i=1}^k\|\p_x^{k+1-i}w\|^2_\infty\int\x^n|\p_x^iw_\eta|^2,\label{2.22}\\
&\int\x^n\sum_{i=1}^k|\p_x^{k+1-i}w||\p_x^iw_\eta||\p_x^{k+1}w|\nonumber\\
&\leq C\int\x^n|\p_x^{k+1}w|^2+C\sum_{i=1}^k\|\p_x^{k+1-i}w\|^2_\infty\int\x^n|\p_x^iw_\eta|^2.\label{2.23}
\end{align}
For case $i=0$, there holds
\begin{align}
&\int\x^n|\p_x^{k+1}w||w_\eta||\p_x^{k+1}w_\eta|\leq \delta\int\x^n|\p_x^{k+1}w_\eta|^2+C\|w_\eta\|^2_\infty\int\x^n|\p_x^{k+1}w|^2,\label{2.24}\\
&\int\x^n|\p_x^{k+1}w|^2|w_\eta|\leq C(1+|w_\eta\|^2_\infty)\int\x^n|\p_x^{k+1}w|^2.\label{2.25}
\end{align}
Substituting (\ref{2.21})-(\ref{2.25}) into (\ref{2.20}), together with using (\ref{2.15}), we obtain
\begin{align}\label{2.26}
&\frac{d}{dx}\int\x^n|\p_x^{k+1}w|^2+\int\x^n|\p_x^{k+1}w_\eta|^2\nonumber\\
&\lesssim (1+\mathcal{N}_1(x))^2\int\x^n|\p_x^{k+1}w|^2+\sum_{i=1}^k\|\p_x^{k+1-i}w\|^2_\infty\int\x^n|\p_x^iw_\eta|^2.
\end{align}
It remains to give bound on $\|\p_x^iw\|^2_\infty$ for $1\leq i\leq k$. Recalling that $w$ vanishes on $\eta=0$ and $\eta=\infty$. Then there holds
\begin{align*}
|\p_x^iw|^2&=\int_0^\eta\p_\eta(|\p_x^iw|^2)\leq \int|\p_x^iw||\p_x^iw_\eta|
\leq \int|\p_x^iw|^2+\int|\p_x^iw_\eta|^2,\nonumber\\
&\leq \int|\p_x^iw|^2+\int|\p_x^iw_\eta(0,\eta)|^2+\int_0^x\p_x\int|\p_x^iw_\eta|^2\nonumber\\
&\leq \mathcal{N}_{i+1}(x)+\int|\p_x^iw_\eta(0,\eta)|^2,
\end{align*}
which gives
\begin{align}\label{2.27}
\|\p_x^iw\|^2_\infty\leq \mathcal{N}_{i+1}(x)+\int|\p_x^iw_\eta(0,\eta)|^2.
\end{align}
For the estimate to $\p_x^iw_\eta(0,\eta)$, we should also prove by mathematical induction. Indeed, for $i=1$, multiplying (\ref{2.11}) by $w_x$ and integrating by parts over $\mathbb{R}_+$ gives
\begin{align}\label{2.28}
\int|\p_xw_\eta|^2\lesssim \int[|w_{xx}||w_x|+|\p_xw_\eta||w_x||w_\eta|+|w_x||\p_xw_\eta|],
\end{align}
applying Cauchy's inequality to which implies that
\begin{align}\label{2.29}
\int|\p_xw_\eta|^2\lesssim (1+\|w_\eta\|^2_\infty)\int|w_x|^2+\int|w_{xx}|^2\leq C(\mathcal{N}_2(x)+1)^2
\end{align}
Taking $x\rightarrow 0$ yields $\int|\p_xw_\eta(0,\eta)|^2\leq C(\mathcal{N}_2(0)+1)^2$. Next, assume that there holds
\begin{align}\label{2.30}
\sum^{i-1}_{\alpha=1}\int|\p^\alpha_xw_\eta(0,\eta)|^2\leq C(\mathcal{N}_i(0)+1)^2
\end{align}
for $i\geq2$. Then, similarly applying $\p_x^i$ to (\ref{2.6})$_1$, multiplying the result by $\p_x^iw$ and integrating by parts over $\mathbb{R}_+$, we have
\begin{align}\label{2.31}
&\int|\p_x^iw_\eta|^2\lesssim (1+\|w_\eta\|^2_\infty)\mathcal{N}_{i+1}(x)+\sum_{\alpha=1}^{i-1}\|\p_x^{i-\alpha}w\|^2_\infty\int|\p_x^\alpha w_\eta|^2\nonumber\\
&\lesssim (\mathcal{N}_{i+1}(x)+1)^2+\sum_{\alpha=1}^{i-1}\int(|\p_x^{i-\alpha}w|^2+|\p_x^{i-\alpha}w_\eta|^2)\int|\p_x^\alpha w_\eta|^2.
\end{align}
Hence, taking $x\rightarrow 0$ in (\ref{2.31}), together with (\ref{2.30}), we get
\begin{align}\label{2.32}
\sum^i_{\alpha=1}\int|\p^\alpha_xw_\eta(0,\eta)|^2\leq C(\mathcal{N}_{i+1}(0)+1)^2.
\end{align}
Therefore, by mathematical induction, (\ref{2.32}) holds for any $i\geq 1$. Now, substituting (\ref{2.32}) into (\ref{2.27}) and further substituting (\ref{2.27}) into (\ref{2.26}), we have
\begin{align}\label{2.33}
\frac{d}{dx}\int\x^n|\p_x^{k+1}w|^2+\int\x^n|\p_x^{k+1}w_\eta|^2\lesssim \mathcal{N}_{k+1}(x)+(\mathcal{N}_{k+1}(x)+1)\sum_{i=1}^k\int\x^n|\p_x^iw_\eta|^2.
\end{align}
Finally, integrating (\ref{2.33}) over $[0,x]$, add the result to (\ref{2.18}) and using Gronwall's inequality give
\begin{align}\label{2.34}
\mathcal{N}_{k+1}(x)\leq C(\mathcal{N}_{k+1}(0)+1),
\end{align}
and hence, by mathematical induction, (\ref{2.18}) is valid for any $k\geq 0.$

Basing on the solvability of system (\ref{2.6}) and the estimates (\ref{2.34}) for the solution, we are able to prove the solvability of (\ref{2.4}) and the estimates for solution $u^\infty_p$. Precisely, we prove the following:

\begin{prop}\label{p2.2}
Assume that $u^\infty_p(0,y):=\bar{u}_0(y)$ is smooth. Then there exists an smooth solution $u^\infty_p$ to system (\ref{2.4}) satisfies that for any $n,k\in\mathbb{N}$
\begin{align}\label{2.0}
\sup_{x\in[0,L]}\|\y^n\p_x^ku_p^\infty\|_{L^2(\mathbb{R}_+)}+\|\y^n\p_x^ku^\infty_{py}\|_{L^2(0,L;L^2(\mathbb{R}_+))}\leq C_0(n,k,\bar{u}_0).
\end{align}
\end{prop}

\Proof In view of the definition of $w$, we obtain that there exists an unique solution $u_p^\infty(x,y(\eta))$ $=w(x,\eta)+[u_b-u_e]e^{-\eta}$ satisfies (\ref{2.4}) on $[0,L]\times\mathbb{R}_+$. Moreover, since that $u_e+u_p^\infty$ is positively bounded from lower and upper, $\eta$ is equivalent to $y$. Therefore (\ref{2.0}) follows from (\ref{2.34}) and  the reversibility of the von Mises transformation.
\endProof

\begin{cor}\label{c2.3}
Let $u^\infty_p$ be constructed in Proposition \ref{p2.2}, and $v_p^\infty$ be obtained directly by the divergence-free condition. Then, there holds
\begin{align}\label{2.00}
\sup_{x\in[0,L]}\|\y^n\p_x^k\p_y^j[u_p^\infty,v_p^\infty]\|_{L^2(\mathbb{R}_+)}\leq C_0(n,k,j,\bar{u}_0),
\end{align}
for any given $n,j,k\in\mathbb{N}$.
\end{cor}

\Proof
Clearly, (\ref{2.0}) gives the estimate of $u_p^\infty$ in (\ref{2.00}) with $j=0$.

Applying $\p_y$ to equation (\ref{2.4})$_1$ implies that $u^\infty_{py}$  satisfies
\begin{align}\label{2.01}
(u_e+u_p^\infty)u^\infty_{pyx}-\int_0^yu^\infty_{px}d\theta u^\infty_{pyy}=u^\infty_{pyyy}.
\end{align}
In addition, in view of (\ref{2.4})$_2$, we obtain the following boundary conditions
\begin{align}\label{2.02}
u^\infty_{pyy}(x,0)=0,\ \lim_{y\rightarrow 0}u^\infty_{py}(x,y)=0, u^\infty_{py}(0,y)=\bar{u}^\prime_{0}(y).
\end{align}
Then, applying $\p_x^k$ to (\ref{2.01}), multiplying the result by $\p_x^ku^\infty_{py}\y^{2n}$ and integrating by parts over $\mathbb{R}_+$ yield
\begin{align}\label{2.03}
&\frac{1}{2}\frac{d}{dx}\int(u^\infty_p+u_e)|\p_x^ku^\infty_{py}|^2\y^{2n}+\int|\p_x^ku^\infty_{pyy}|^2\y^{2n}\nonumber\\
=&-\frac{1}{2}\int\int_0^yu^\infty_{px}d\theta|\p_x^ku^\infty_{py}|^2\p_y[\y^{2n}]
-\int\p_x^ku^\infty_{pyy}\p_x^ku^\infty_{py}\p_y[\y^{2n}]\nonumber\\
&-\sum_{\ell=0}^{k-1}C_k^\ell\int\left[\p_x^{k-\ell}u^\infty_p\p_x^\ell u_{pxy}^\infty+\int_0^y\p_x^{k-\ell}u^\infty_{px}d\theta\p_x^\ell u_{pyy}^\infty\right]\p_x^ku^\infty_{py}\y^{2n}\nonumber\\
=&:\mathcal{I}_1+\mathcal{I}_2+\mathcal{I}_3+\mathcal{I}_4.
\end{align}

Note that
\begin{align}
\mathcal{I}_1+\mathcal{I}_2&\leq \frac{1}{4}\|\p_x^ku^\infty_{pyy}\y^n\|^2_2 +C\left(\|u^\infty_{px}\y^n\|_2+1\right)\|\p_x^ku^\infty_{py}\y^n\|^2_2,\label{2.04}\\
\mathcal{I}_3&\leq C\sum_{\ell=0}^{k-1}\|\p_x^{k-\ell}u^\infty_{py}\y^n\|_2\|\p_x^{\ell+1} u^\infty_{py}\y^n\|_2\|\p_x^ku^\infty_{py}\y^n\|_2\nonumber\\
&\leq C\sum_{\ell=1}^{k}\|\p_x^{\ell}u^\infty_{py}\y^n\|^2_2+C\sum_{\ell=1}^k\|\p_x^{\ell}u^\infty_{py}\y^n\|
^2_2\|\p_x^ku^\infty_{py}\y^n\|^2_2,\label{2.06}\\
\mathcal{I}_4&=\sum_{\ell=0}^{k-1}\int\p_x^{k-\ell}u^\infty_{px}\p_x^\ell u_{py}^\infty\p_x^ku^\infty_{py}\y^{2n}
+\sum_{\ell=0}^{k-1}\int\int_0^y\p_x^{k-\ell}u^\infty_{px}d\theta\p_x^\ell u^\infty_{py}\p_y\left[\p_x^ku^\infty_{py}\y^{2n}\right]\nonumber\\
&\leq C\sum_{\ell=0}^{k-1}\|\p_x^{k+1-\ell}u^\infty_{p}\y^n\|_2\|\p_x^{\ell}u^\infty_{py}\y^n\|_2
\|\p_x^ku^\infty_{pyy}\y^n\|_2\nonumber\\
&\quad+C\sum_{\ell=0}^{k-1}\|\p_x^{k+1-\ell}u^\infty_p\y^n\|_2\|\p_x^{\ell}u^\infty_{py}\y^n\|_2
\|\p_x^ku^\infty_{py}\y^n\|_2\nonumber\\
&\leq \frac{1}{4}\|\p_x^ku^\infty_{pyy}\y^n\|^2_2+C\|\p_x^ku^\infty_{py}\y^n\|^2_2
+C\sum_{\ell=0}^{k-1}\|\p_x^{k+1-\ell}u^\infty_p\y^n\|^2_2\|\p_x^{\ell}u^\infty_{py}\y^n\|^2_2.\label{2.07}
\end{align}
Substituting (\ref{2.04})-(\ref{2.07}) into (\ref{2.03}), applying Gronwall's inequality and using (\ref{2.0}) and the positivity condition $u_e+u^\infty_p\geq c_0$, we have
\begin{align}\label{2.08}
\sup_{x\in[0,L]}\|\y^n\p_x^ku^\infty_{py}\|_{L^2(\mathbb{R}_+)}+\|\y^n\p_x^ku^\infty_{pyy}\|_{L^2(0,L;L^2(\mathbb{R}_+))}\leq C(n,k,\bar{u}_0).
\end{align}
This gives the estimate of $u_p^\infty$ in (\ref{2.00}) with $j=1$.

Similarly, applying $\p_x^k$ to (\ref{2.3}) yields
\[\p_x^ku^\infty_{pyy}=\sum_{\ell=0}^kC_k^\ell\p_x^{k-\ell}(u_e+u^\infty_p)\p_x^\ell u^\infty_{px}
+\sum_{\ell=0}^kC_k^\ell\int_0^y\p_x^{k-\ell}u^\infty_{px}d\theta \p_x^\ell u^\infty_{py}.\]
Direct calculation gives the estimate of $u_p^\infty$ in (\ref{2.00}) with $j=2$, where (\ref{2.08}) has been used. Then, by iteration method, the estimate of $u_p^\infty$ in (\ref{2.00}) can be derived with arbitrary $j$.

With the estimates of $u_p^\infty$ in hand, we are able to derive the estimates for $v^\infty_p$. In view of the divergence-free condition, we have
\[|\p_x^kv^\infty_p|^2=\left|\int_y^\infty\p_x^{k+1}u^\infty_{p}d\theta\right|^2\leq C\|\p_x^{k+1}u^\infty_{p}\y^m\|^2_2\y^{2-2m}, \]
for any $m\in\mathbb{N}$.
This together with (\ref{2.0}) implies that
\[\|\p_x^kv^\infty_p\y^n\|^2_2\leq C\|\p_x^{k+1}u^\infty_{p}\y^m\|^2_2\int\y^{2n+2-2m}\leq C(k,n,\bar{u}_0),\]
where we take $m=n+2.$

Finally, for any $j\geq 1$, since that $\p_x^k\p_y^jv^\infty_p=-\p_x^{k+1}\p_y^{j-1}u^\infty_p$, the proof of (\ref{2.00}) is completed directly by the established estimates of $u_p^\infty.$
\endProof

\begin{prop}\label{p2.4}
Under the assumptions in Theorem 1.1, there exists smooth functions $[u^0_p,v^0_p]$, defined in $\Omega_\e$, satisfying the following inhomogeneous system:
 \begin{align}\label{2.09}
 \begin{cases}
 (u_e+u_p^0)u^0_{px}+(v^0_p+v_e^1(x,0))u^0_{py}-u^0_{pyy}=R^{u,0}_p,\\
u^0_{px}+v^0_{py}=0,\\
u^0_p(x,0)=u_b-u_e,\ u^0_{py}(x,\frac{1}{\sqrt\e})=0,\ [v_p^0+v_e^1](x,0)=0,\ v^0_p(x,\frac{1}{\sqrt\e})=0,
 \end{cases}
 \end{align}
 where the inhomogeneous term $R^{u,0}_p$ is a higher order term of $\sqrt\e$. In addition, it holds that
\begin{align}\label{2.010}
\sup_{x\in[0,L]}\|\y^n\p_x^k\p_y^j[u_p^0,v_p^0]\|_{L^2(I_\e)}\leq C_0(n,k,j,\bar{u}_0),
\end{align}
for any given $n,j,k\in\mathbb{N}$.
\end{prop}
\Proof Let $u^\infty_p$ be constructed in Proposition \ref{p2.2}, and $v_p^\infty$ be obtained directly by the divergence-free condition. Define that
\begin{align}\label{2.011}
&u^0_p(x,y):=\chi(\sqrt\e y)u^\infty_p(x,y)-\sqrt\e\chi^\prime(\sqrt\e y)\int^\infty_yu^\infty_p(x,\theta)d\theta,\nonumber\\
&v^0_p(x,y):=\chi(\sqrt\e y)v^\infty_p(x,y).
\end{align}
Then, it follows from directly calculation that $[u^0_p,v^0_p]$ satisfies (\ref{2.09}) with
\begin{align}\label{2.012}
R^{u,0}_p=&\sqrt\e\chi\int_0^y\chi^\prime d\theta(u^\infty_pu^\infty_{px}+v^\infty_pu^\infty_{py})
-\sqrt\e\chi^\prime\chi u^\infty_{px}\int^\infty_y u^\infty_p d\theta\nonumber\\
&-\sqrt\e\chi^\prime v^\infty_p(u_e+\chi u^\infty_p)-3\sqrt\e\chi^\prime u^\infty_{py}+2\sqrt\e\chi^\prime u^\infty_p\int_0^y\chi v^\infty_{py}d\theta\nonumber\\
&+2\e\chi^\prime u_p^\infty\int_0^y\chi^\prime v^\infty_pd\theta-3\e\chi^{\prime\prime}u^\infty_p+\e(\chi^\prime)^2v^\infty_p\int_y^\infty u^\infty_pd\theta\nonumber\\
&-\e\chi^{\prime\prime}(\chi v^\infty_p-v_p^\infty(0))\int_y^\infty u^\infty_pd\theta+\e^{3/2}\chi^{\prime\prime\prime}\int_y^\infty u^\infty_p d\theta\nonumber\\
:=&\sqrt\e E_1+\e E_2.
\end{align}
Finally, using (\ref{2.00}) together with the definition of $\chi(\cdot)$ give estimate (\ref{2.010}).
\endProof

\section{The first order Euler corrector}

To construct the first order Euler corrector $[u^1_e,v^1_e,p_e^1]$, we first formulate a closed system for these functions. For one hand, denote
\begin{align*}
R_1^u:=&(u_e^1+u_p^1)u_{px}^0+(u^0_e+u_p^0)(u_{ex}^1+u_{px}^1)+v_p^1\p_y(u_e^0+u_p^0)+(v_p^0+v_e^1)\p_y(u_e^1+u_p^1)\\
&+(p^1_{ex}+p^1_{px})-\p_y^2(u_e^1+u^1_p)+(v_p^0+v^1_e)u^0_{ez}+y(u^0_{ez}u^0_{px}+v^1_{ez}u^0_{py})+E_1.
\end{align*}
Note that
\[(v^0_p+v^1_e)\p_yu^1_e=\sqrt\e(v^0_p+v^1_e)u^1_{ez},\ \p_y^2u^1_e=\e u^1_{ezz}.\]
Since the unknown Euler corrector $[u_e^1,v_e^1,p_e^1]$ and Prandtl corrector $[u_p^1,v_p^1,p_p^1]$ couple with each other, we take equation
\begin{align}\label{3.1}
u^0_eu^1_{ex}+v^1_eu^0_{ez}+p^1_{ex}=0,
\end{align}
for the first order Euler corrector and when it has been constructed, we take
\begin{align}\label{3.2}
&(u_e^1+u_p^1)u^0_{px}+(u^0_e+u^0_p)u^1_{px}+u_p^0u^1_{ex}+(v^1_e+v^0_p)u^1_{py}+v^1_p\p_y(u^0_e+u^0_p)\nonumber\\
&+p^1_{px}-u^1_{pyy}+v^0_pu^0_{ez}+y(u_{ez}^0u_{px}^0+v^1_{ez}u^0_{py})+E_1=0.
\end{align}
for the first order Prandtl corrector. Hence, the error $R^u_1$ then reads
\begin{align}\label{3.201}
\sqrt\e(v^0_p+v^1_e)u^1_{ez}-\e u^1_{ezz}.
\end{align}

On the other hand, in view of the divergence-free condition, we have
\begin{align}
&u^1_{ex}+v^1_{ez}=0,\label{3.3}\\
&u^1_{px}+v^1_{py}=0.\label{3.4}
\end{align}
Even so, the equations above are still not enough to construct neither $[u_e^1,v_e^1,p_e^1]$ or $[u_p^1,v_p^1,p_p^1]$. This motivates us to consider the vertical component (\ref{1.7}).

Denote that
\begin{align*}
R_0^v:=(u_e^0+u_p^0)(v^0_{px}+v_{ex}^1)+(v_p^0+v_e^1)\p_y(v_p^0+v_e^1)+p^1_{ez}+\frac{p^1_{py}}{\sqrt\e}+p^2_{py}-\p_y^2(v_p^0+v^1_e)
\end{align*}
Clearly, the leading term in $R_0^v$ is $p^1_{py}$. Let $p^1_{py}=0$, that is,
\begin{align}\label{3.5}
p_p^1=p_p^1(x).
\end{align}
Similar to (\ref{3.1}) and (\ref{3.2}), we take
\begin{align}\label{3.6}
u_e^0v_{ex}^1+p^1_{ez}=0,
\end{align}
and
\begin{align}\label{3.7}
(u_e^0+u_p^0)v^0_{px}+u_p^0v_{ex}^1+(v_p^0+v_e^1)v^0_{py}+p^2_{py}-v^0_{pyy}=0.
\end{align}
Then the error $R^v_0$ is reduced to
\begin{align}\label{3.701}
\sqrt\e(v_p^0+v^1_e)v^1_{ez}-\e v^1_{ezz}.
\end{align}

In conclusion, we get a system consisting of (\ref{3.1}), (\ref{3.3}) and (\ref{3.6}) to construct $[u_e^1,v_e^1,p_e^1]$, and another system consisting of (\ref{3.2}), (\ref{3.4}) and (\ref{3.5}) to construct $[u_p^1,v_p^1,p_p^1]$. After these functions being given, $p^2_p$ will be determined directly by (\ref{3.7}).

In this section, we only focus on the construction of $[u_e^1,v_e^1,p_e^1]$, while the construction of $[u_p^1,v_p^1,p_p^1]$ will be done in the next section.

Eliminating $p^1_e$ in (\ref{3.1}) and (\ref{3.6}) by using equation (\ref{3.3}), we deduce the following elliptic equation for $v^1_e$:
\begin{align}\label{3.8}
-u_e^0\Delta v^1_e+u^0_{ezz}v^1_e=0,\ \textrm{in}\ \Omega_0,
\end{align}
with $\Delta:=\p_x^2+\p_z^2.$ In order to solve this equation, we take the following boundary conditions
\begin{align}\label{3.9}
v_e^1(x,0)=-v^0_p(x,0),\ v^1_e(x,1)=0,\ v^1_e(0,z)=V_{b0}(z),\ v^1_e(L,z)=V_{bL}(z),
\end{align}
with the compatibility assumption:
\begin{align}\label{3.10}
 [V_{b0}(0),V_{bL}(0)]=-[v_p^0(0,0),v_p^0(L,0)]\ \mathrm{and}\ V_{b0}(1)=V_{bL}(1)=0.
\end{align}

To avoid singularity caused by the presence of corners in $\Omega_0$, we instead consider the modified elliptic problem:
\begin{align}\label{3.11}
-u_e^0\Delta v^1_e+u^0_{ezz}v^1_e=E_b, \ \textrm{in}\ \Omega_0,
\end{align}
with boundary conditions (\ref{3.9}). Later, we shall construct a proper potential $E_b$ such that $v^1_e$, the solution of the elliptic equation (\ref{3.11}), is regular enough and that $\int_z^\infty E_bd\theta\rightarrow 0$ as $\e\rightarrow 0$.

To define $E_b$, we first introduce
\begin{align}\label{3.12}
B(x,z):=\left(1-\frac{x}{L}\right)\frac{V_{b0}(z)}{v^0_p(0,0)}v^0_p(x,0)+\frac{x}{L}\frac{V_{bL}(z)}{v^0_p(L,0)}v^0_p(x,0),
\end{align}
in the case of both $v^0_p(0,0)$ and $v^0_p(L,0)$ are nonzero, while in the case that $v^0_p(0,0)=0$ or $v^0_p(L,0)=0$, we simply replace the ratio $\frac{V_{b0}(z)}{v^0_p(0,0)}v^0_p(x,0)$ by $V_{b0}(z)-v^0_p(x,0)(1-z)$, or replace $\frac{V_{bL}(z)}{v^0_p(L,0)}v^0_p(x,0)$ by $V_{bL}(z)-v^0_p(x,0)(1-z)$, respectively. We infer that $B(x,z)$ satisfies all the boundary conditions in (\ref{3.9}).

Then denote $F_e(x,z):=-u_e^0\Delta B+u^0_{ezz}B.$ In view of the estimates (\ref{2.00}), it is clear that $B\in W^{k,p}(\Omega_0)$ for arbitrary $k\geq0,p> 1$, provided $V_{b0}(z),V_{bL}(z)\in W^{k,p}(0,1)$, and hence $F_e\in W^{k,p}(\Omega_0)$.

Now, take $E_b=:\chi(\frac{z}{\e})F_e(x,0)$. Before solving equation (\ref{3.11}) and derive estimates for $v^1_e$, we consider the following auxiliary problem
\begin{align}\label{3.13}
\begin{cases}
-u^0_e\Delta \tilde{w}+u^0_{ezz}\tilde{w}=E_b-F_e, \textrm{in}\ \Omega_0,\\
\tilde{w}\big|_{\p\Omega_0}=0.
\end{cases}
\end{align}
Precisely, we prove the following lemma:

\begin{lem}\label{L3.1}
Assume that $F_e(x,z)\in W^{k,p}(\Omega_0)$ for any $k\geq 0,p>1$. Then there exists a unique smooth solution to the boundary value problem (\ref{3.13}) satisfying that
\begin{align}\label{3.131}
\|\tilde{w}\|_{L^\infty(\Omega_0)}+\|\tilde{w}\|_{H^2(\Omega_0)}\leq C,\ \ \|\tilde{w}\|_{H^3(\Omega_0)}\leq C\e^{-\frac{1}{2}},\ \ \|\tilde{w}\|_{H^4(\Omega_0)}\leq C\e^{-\frac{3}{2}},
\end{align}
where $C$ is independent of $\e$. In addition,there holds
\begin{align}\label{3.132}
\|\tilde{w}\|_{W^{2,q}(\Omega_0)}\leq C,\ \ \|\tilde{w}\|_{W^{3,q}(\Omega_0)}\leq C\e^{-1+\frac{1}{q}},\ \ \|\tilde{w}\|_{W^{4,q}(\Omega_0)}\leq C\e^{-2+\frac{1}{q}},
\end{align}
for any $q\in(1,+\infty)$ and $C$ being independent of $\e$.
\end{lem}

\Proof
Define bilinear form on $H^1_0(\Omega_0)$:
\[\mathcal{B}[\tilde{w},\tilde{v}]:=\iint_{\Omega_0}\left(\nabla\tilde{w}\cdot\nabla\tilde{v}
+\frac{u^0_{ezz}}{u^0_e}\tilde{w}\tilde{v}\right).\]
Note that, on one hand,
\begin{align*}
\int_0^1|\p_z\tilde{w}|^2&=\int_0^1\left|\p_z\left(\frac{\tilde{w}}{u_e^0}u^0_e\right)\right|^2
=\int_0^1\left|\p_z\left(\frac{\tilde{w}}{u_e^0}\right)u^0_e+\frac{\tilde{w}}{u^0_e}u^0_{ez}\right|^2\\
&=\int_0^1\left|\p_z\left(\frac{\tilde{w}}{u_e^0}\right)\right|^2|u^0_e|^2+\int_0^1\left|\frac{\tilde{w}}{u^0_e}\right|^2 |u^0_{ez}|^2+2\int_0^1\p_z\left(\frac{\tilde{w}}{u_e^0}\right)\frac{\tilde{w}}{u^0_e}u^0_e u^0_{ez}\\
&=\int_0^1\left|\p_z\left(\frac{\tilde{w}}{u_e^0}\right)\right|^2|u^0_e|^2+\int_0^1\left|\frac{\tilde{w}}{u^0_e}\right|^2 |u^0_{ez}|^2-\int_0^1\left|\frac{\tilde{w}}{u_e^0}\right|^2[u^0_e u^0_{ez}]_z\\
&=\int_0^1\left|\p_z\left(\frac{\tilde{w}}{u_e^0}\right)\right|^2|u^0_e|^2-\int_0^1\frac{u^0_{ezz}}{u^0_e}\tilde{w}^2,
\end{align*}
which implies that
\begin{align}\label{3.14}
\iint_{\Omega_0}\left(|\p_z\tilde{w}|^2+\frac{u^0_{ezz}}{u^0_e}\tilde{w}^2\right)
=\iint_{\Omega_0}\left|\p_z\left(\frac{\tilde{w}}{u_e^0}\right)\right|^2|u^0_e|^2.
\end{align}
On the other hand, thanks to the positivity and smoothness of $u_e^0$, we have
\begin{align}\label{3.15}
\iint_{\Omega_0}|\p_z\tilde{w}|^2&=\iint_{\Omega_0}\left|\p_z\left(\frac{\tilde{w}}{u_e^0}u^0_e\right)\right|^2
\leq 2\iint_{\Omega_0}\left|\p_z\left(\frac{\tilde{w}}{u_e^0}\right)\right|^2|u^0_e|^2
+2\iint_{\Omega_0}\left|\frac{\tilde{w}}{u^0_e}\right|^2|u^0_{ez}|^2\nonumber\\
&\leq 2\iint_{\Omega_0}\left|\p_z\left(\frac{\tilde{w}}{u_e^0}\right)\right|^2|u^0_e|^2
+2\int_0^L\left[\int_0^1\left|\p_z\left(\frac{\tilde{w}}{u^0_e}\right)\right|^2\int_0^1|u^0_{ez}|^2\right]\nonumber\\
&\leq 2\iint_{\Omega_0}\left|\p_z\left(\frac{\tilde{w}}{u_e^0}\right)\right|^2|u^0_e|^2
+2\|u^0_{ez}\|_{L^2(0,1)}^2\iint_{\Omega_0}\left|\p_z\left(\frac{\tilde{w}}{u^0_e}\right)\right|^2\nonumber\\
&\leq 2\iint_{\Omega_0}\left|\p_z\left(\frac{\tilde{w}}{u_e^0}\right)\right|^2|u^0_e|^2
\left[1+\frac{\|u^0_{ez}\|_{L^2(0,1)}^2}{\min_z|u^0_e|^2}\right]\nonumber\\
&\leq C_0\iint_{\Omega_0}\left|\p_z\left(\frac{\tilde{w}}{u_e^0}\right)\right|^2|u^0_e|^2,
\end{align}
where the inequality $|f(z)|\leq \sqrt z\|f_z(z)\|_2 $ has been used. Then, by applying the Poincar\'e inequality, it follows that
\begin{align}\label{3.16}
\mathcal{B}[\tilde{w},\tilde{w}]\geq \alpha\|\tilde{w}\|^2_{H^1},
\end{align}
where $\alpha$ is a positive constant. In addition, by applying the Cauchy inequality, there holds that
\begin{align}\label{3.17}
\mathcal{B}[\tilde{w},\tilde{v}]\leq\beta\|\tilde{w}\|_{H^1}\|\tilde{v}\|_{H^1},
\end{align}
for any $\tilde{w},\tilde{v}\in H_0^1(\Omega_0)$.
Moreover, since $F_e\in W^{k,p}(\Omega_0),\ k\geq 0,\ p\geq 1$, it is clear that
 \[\|\p_z^kE_b\|_p\leq C\e^{-k+\frac{1}{p}}.\]
 Therefore, by Lax-Milgram theorem, there exists an unique weak solution $\tilde{w}\in H^1_0(\Omega_0)$ to problem (\ref{3.13}) satisfying $\|\tilde{w}\|_{H^1}\leq C.$

Now, rewrite (\ref{3.13}) as below:
\begin{align}\label{3.19}
\begin{cases}
-\Delta \tilde{w}=G_e=:\left(E_b-F_e-u^0_{ezz}\tilde{w}\right)/u^0_e, \textrm{in}\ \Omega_0,\\
\tilde{w}\big|_{\p\Omega_0}=0.
\end{cases}
\end{align}
Clearly, $G_e\in L^2(\Omega_0)$. Then by the elliptic estimates, we have
 \[\|\tilde{w}\|_{H^2}\leq C\|G_e\|_2\leq C.\]
 In addition, since $\tilde{w}=0$ on the boundary, we obtain that
\begin{align*}
&|\tilde{w}(x,z)|\leq \int_0^x|\tilde{w}_x(s,z)|ds\leq 2\int_0^x\left(\int_0^z|\tilde{w}_x\tilde{w}_{xz}|d\theta\right)^{1/2}ds\\
&\leq 2\sqrt x\|\tilde{w}_x\|^{1/2}_{L^2(\Omega_0)}\|\tilde{w}_{xz}\|^{1/2}_{L^2(\Omega_0)}\leq C\sqrt L,
\end{align*}
which implies the uniform boundedness of $\tilde{w}$.

Next, we derive the higher regularity estimates for $\tilde{w}$. Since that $E(x,0)-F(x,0)=0$, we get $G_e(x,0)=0$ and hence, by equation (\ref{3.19}), $\tilde{w}=\tilde{w}_{zz}=0$ on $z=0$. Then, applying $\p_z$ and $\p_{zz}$ to (\ref{3.19}) yields the elliptic problems for $\tilde{w}_z$ and $\tilde{w}_{zz}$, respectively:
\begin{align}\label{3.20}
\begin{cases}
-\Delta \tilde{w}_z=\p_z\left[(E_b-F_e-u^0_{ezz}\tilde{w})/u^0_e\right], \textrm{in}\ \Omega_0,\\
\tilde{w}_z\big|_{x=0,L}=\tilde{w}_{zz}\big|_{z=0}=0,\ \tilde{w}_{zz}\big|_{z=1}=0,
\end{cases}
\end{align}
and
\begin{align}\label{3.21}
\begin{cases}
-\Delta \tilde{w}_{zz}=\p^2_z\left[(E_b-F_e-u^0_{ezz}\tilde{w})/u^0_e\right], \textrm{in}\ \Omega_0,\\
\tilde{w}_{zz}\big|_{x=0,L}=\tilde{w}_{zz}\big|_{z=0}=0,\ \tilde{w}_{zz}\big|_{z=1}=0,
\end{cases}
\end{align}
where the higher order compatibility condition $V_{b0}^{\prime\prime}(1)=V^{\prime\prime}_{bL}(1)=0$ has been used.

Again, the elliptic estimates for $H^2$ norm and the estimates for $E_b$ and $F_e$ then give
\[\|\p_z^k\tilde{w}\|_{H^2}\leq C\e^{-k+\frac{1}{2}},\ k=1,2.\]
To complete the $H^3$ and $H^4$ estimates for $\tilde{w}$, it remains to estimate $L^2$ and $H^1$ norm for $\tilde{w}_{xxx}$. Applying $\p_x$ to equation (\ref{3.19}), we have
\[-\tilde{w}_{xxx}=\tilde{w}_{zzx}+\p_x\left[(E_b-F_e-u^0_{ezz}\tilde{w})/u^0_e\right],\]
which give the $L^2$ and $H^1$ norm estimates for $\tilde{w}_{xxx}$ and hence completes the proof of (\ref{3.131}).

Finally, the $W^{k,q}$ estimates follow simply from the standard elliptic theory. The proof of this lemma is completed.
\endProof

Now, take $v^1_e=B+\tilde{w}$. Then, recalling that $B$ satisfies boundary conditions (\ref{3.9}), it follows that $v^1_e\in W^{k,p}(\Omega_0)$ is the unique smooth solution to equation (\ref{3.11}) with boundary conditions (\ref{3.9}). It should be noted that $v^1_{ezz}(x,1)=0$, since the definition of $E_b$ and the equation (\ref{3.11}). In addition, as $B\in W^{4,q}(\Omega_0)$, there holds
\begin{align}\label{3.22}
\|v^1_e\|_\infty+\|v^1_e\|_{W^{2,q}}\leq C,\ \|v^1_e\|_{W^{2+k,q}}\leq C\e^{-k+\frac{1}{q}},k=1,2.
\end{align}

Furthermore, in view of equation (\ref{3.3}) and (\ref{3.6}), we take
\begin{align*}
&u^1_e(x,z)=u^1_b(z)-\int_0^xv^1_{ez}(\xi,z)d\xi,\\
&p^1_e(x,y)=\int_z^1 u^0_e(\theta)v^1_{ex}(x,\theta)d\theta-\int_0^xu^0_e(1)v^1_{ez}(s,1)ds,
\end{align*}
where $u^1_b(z)=u^1_e(0,z)$ satisfies $\p_zu^1_b(1)=0$, and hence we have $u^1_{ez}(x,1)=0$.

Substituting $u^1_e,p^1_e$ into (\ref{3.1}) and integrating by parts yield
\[u^0_eu^1_{ex}+u^0_{ez}v_e^1+p^1_{ex}=\int_z^1(u^0_e\Delta v^1_e-u^0_{ezz}v_e^1)d\theta=-\int_z^1 E_b(x,\theta)d\theta.\]

Base on the estimates for $v^1_e$ and $E_b$, we infer that
\begin{align*}
&\|u^1_{e}\|_{L^\infty(\Omega_\e)}\leq C,\ \ \|u^1_e\|_{H^1(\Omega_\e)}\leq C\e^{-\frac{1}{4}},\\
&\|u^1_{ezz}\|_{L^2(\Omega_\e)}\leq C\e^{-\frac{1}{4}}\|u^1_{bzz}\|_{L^2(0,1)}+C\e^{-\frac{1}{4}}\|v^1_{ezzz}\|_{L^2(\Omega_0)}\leq C\e^{-\frac{3}{4}}\\
&\|(v_p^0+v_e^1)u^1_{ez}\|_{L^2(\Omega_\e)}\leq\|v_p^0+v_e^1\|_{L^\infty(\Omega_\e)}\|u^1_{ez}\|_{L^2(\Omega_\e)}\leq C\e^{-\frac{1}{4}},\\
&\left\|\int_{\sqrt\e y}^1E_b(x,\theta)d\theta\right\|_{L^2(\Omega_\e)}\leq C\e^{-\frac{1}{4}}\|E_b\|_{L^2(\Omega_0)}\leq C\e^{\frac{1}{4}}.
\end{align*}
Hence, it follows that
\begin{align}
\|R^u_1\|_{L^2(\Omega_\e)}&\leq \sqrt\e\|(v_p^0+v_e^1)u^1_{ez}\|_{L^2(\Omega_\e)}+\e\|\p_z^2 u^1_{e}\|_{L^2(\Omega_\e)}+
\left\|\int_{\sqrt\e y}^1E_b(x,\theta)d\theta\right\|_{L^2(\Omega_\e)}\nonumber\\
&\leq C\e^{\frac{1}{4}},\label{3.23}\\
\|R^v_0\|_{L^2(\Omega_\e)}&\leq \sqrt\e\|v_p^0+v_e^1\|_{L^\infty(\Omega_\e)}\|v^1_{ez}\|_{L^2(\Omega_\e)}
+\e\|v^1_{ezz}\|_{L^2((\Omega_\e)}\leq C\e^{\frac{1}{4}},\label{3.24}
\end{align}

Finally, we estimate $E^0$, which is defined in (\ref{2.1}). Note that
\begin{align*}
&\left|u_{px}^0(x,y)\int_0^y\int_y^ru^0_{ezz}(\sqrt\e\tau)d\tau dr\right|\leq C|u^0_{px}(x,y)|\sup_{z\in[0,1]}|u^0_{ezz}(z)|\y^2,\\
&\left|u_{py}^0(x,y)\int_0^y\int_y^rv^1_{ezz}(x,\sqrt\e\tau)d\tau dr\right|\leq C\e^{-\frac{1}{4}}|u^0_{py}(x,y)|\|v^1_{ezz}\|_{L^2(0,1)}\y^2.
\end{align*}
Then, it follows that
\begin{align}\label{3.25}
\|E^0\|_{L^2(\Omega_\e)}&\leq C\e \|u^0_{px}\y^2\|_{L^2(I_\e)}\sup_{z\in[0,1]}|u^0_{ezz}(z)|
+C\e^{\frac{3}{4}}\|u^0_{py}\y^2\|_{L^2(I_\e)}\|v^1_{ezz}\|_{L^2(\Omega_0)}\nonumber\\
&\leq C\e^{\frac{3}{4}},
\end{align}
where the estimate (\ref{2.00}) has been used.

\section{The first order Prandtl corrector}
In this section, we shall construct the first order Pranndtl corrector $[u^1_p,v^1_p,p^1_p]$, which solves (\ref{3.2}), (\ref{3.4}) and (\ref{3.5}). For convenience, we denote $u^0:=u_e+u^0_p$.

It should be noted that
\begin{align*}
(u^0_e-u_e)u^1_{px}=\sqrt\e u^1_{px}\int_0^yu^0_{ez}(\sqrt\e\theta)d\theta,\ v^1_p\p_yu^0_e=\sqrt\e v^1_pu^0_{ez}.
\end{align*}
Then, (\ref{3.2}) can be rewritten as
\begin{align*}
&u^0u^1_{px}+u^0_xu^1_p+u^0_yv^1_p+(v^0_p+v^1_e)u^1_{py}+p^1_{px}-u^1_{pyy}\\
&=-u^0_{px}u^1_e-u^0_pu^1_{ex}-(v^0_p+yu^0_{px})u^0_{ez}-yu^0_{py}v^1_{ez}-E_1:=F_p.
\end{align*}
We infer that, by Section 2, $u^0$ is positively bounded from both lower and upper. In addition, the error terms should be added to \begin{align}\label{4.0}
\tilde{R}^u_1:=\sqrt\e u^1_{px}\int_0^yu^0_{ez}(\sqrt\e\theta)d\theta+\sqrt\e v^1_pu^0_{ez}.
\end{align}

Taking $p^1_p$ to be an absolute constant implies that $p^1_{px}=0.$ Then the system for $[u^1_p,v^1_p]$ can be rewritten as follows:
\begin{align}\label{4.1}
\begin{cases}
u^0u^1_{px}+u^0_xu^1_p+u^0_yv^1_p+(v^0_p+v^1_e)u^1_{py}-u^1_{pyy}=F_p,\\
u^1_{px}+v^1_{py}=0,
\end{cases}
\end{align}
with the boundary conditions
\[
u^1_p(0,y)=\bar{u}_1(y),u^1_p(x,0)=-u^1_e(x,0),u^1_{py}(x,\frac{1}{\sqrt\e})=v^1_p(x,0)=v^1_p(x,\frac{1}{\sqrt\e})=0.
\]

Similar to the situation in Section 2, we first extend the domain from $I_\e$ to $\mathbb{R}_+$ with the boundary condition $u^1_{py}(x,\frac{1}{\sqrt\e})=v^1_p(x,\frac{1}{\sqrt\e})=0$ being replaced by $u_p(x,\infty)=v_p(x,\infty)=0$, and denote the unknown functions in (\ref{4.1}) by $[u_p,v_p]$, for distinction. To the given functions, define that $u^0_e(z)\equiv u^0_e(1)$,~~$u^1_e(x,z)\equiv u^1_e(x,1)$, $v^1_e(x,z)\equiv 0$ in $z\in(1,+\infty)$, and also $\bar{u}_1(y)\equiv \bar{u}_1(\frac{1}{\sqrt\e})$ in $y\in(\frac{1}{\sqrt\e},+\infty)$. For convenience, we still denote them by $u^0_e,u^1_e,v^1_e,\bar{u}_1$.

Then, applying $\p_y$ to (\ref{4.1})$_1$ and using (\ref{4.1})$_2$ yields
\[-u^0v_{pyy}+u^0_{yy}v_p+u^0_{xy}u_p+(v^0_p+v^1_e)u_{pyy}+\sqrt\e v^1_{ez}u_{py}-u_{pyyy}=F_{py},\]
which also can be rewritten as
\begin{align}\label{4.2}
-v_{pyy}+\frac{u^0_{yy}}{u^0}v_p-\left(\frac{u_{py}}{u^0}\right)_{yy}=G_p,
\end{align}
where we denote
\[G_p:=\frac{1}{u^0}\left[F_{py}-u^0_{xy}u_p-(v^0_p+v^1_e)u_{pyy}-\sqrt\e v^1_{ez}u_{py}\right]-2\left(\frac{1}{u^0}\right)_yu_{pyy}-\left(\frac{1}{u^0}\right)_{yy}u_{py}.\]

Furthermore, applying $\p_x$ to (\ref{4.2}), we get
\begin{align}\label{4.3}
-v_{pxyy}+\frac{u^0_{yy}}{u^0}v_{px}+\left(\frac{v_{pyy}}{u^0}\right)_{yy}=G_{px}-\left(\frac{u^0_{yy}}{u^0}\right)_xv_p
-\left[\left(\frac{1}{u^0}\right)_xu_{py}\right]_{yy}.
\end{align}

The proof of solvability of (\ref{4.1}) on $[0,L]\times\mathbb{R}_+$ consists of several steps.

Step 1,  we establish the estimates for the boundary conditions of $v_p$ in term of the given data $\bar{u}_1(y)$. For simplification, we denote $\|\y^nf\|_{H^k}:=\sum_{i=0}^k\|\y^n\p_y^if\|_2$, for any $n,k\geq 0$ and $f\in H^k.$

\begin{lem}\label{L4.1}
Let $[u_p,v_p]$ be smooth solutions of (\ref{4.1}). Then there holds that
\begin{align}
&\|\y^nv_{pyy}(0,\cdot)\|_{L^2(\mathbb{R}_+)}\leq C_0\left(1+\|\y^n\bar{u}_1\|_{H^3(\mathbb{R}_+)}\right),\label{4.4}\\
&\|\y^nv_{pxyy}(0,\cdot)\|_{L^2(\mathbb{R}_+)}\leq C_0\left(1+\|\y^n\bar{u}_1\|_{H^5(\mathbb{R}_+)}+\|u^1_{exx}(0,\cdot)\|_{H^1(\mathbb{R}_+)}\right),\label{4.5}
\end{align}
for any $n\geq 0$ and some constant $C_0=C_0(u^0, v^0_p,u^1_b,V_{b0})$.
\end{lem}

\Proof
Define stream function $\psi(x,y)=-\int_y^\infty u_p(x,\theta)d\theta$. Then $u_p=\psi_y,v_p=-\psi_x$. Further denote $\phi:=u^0\psi_y-u^0_y\psi$. Then (\ref{4.1})$_1$ becomes
\begin{align}\label{4.6}
\phi_x=-u_{xy}^0\psi-(v^0_p+v^1_e)u_{py}+u_{pyy}+F_p.
\end{align}
By the definition of $\phi$, we have
\begin{align}\label{4.7}
\|\y^n\phi_y(0,\cdot)\|_2\lesssim \|\y^nu^0u_{py}(0,\cdot)\|_2+\|\y^nu^0_{yy}\psi(0,\cdot)\|_2\leq C_0\|\y^n\bar{u}_1\|_{H^1}.
\end{align}
In view of (\ref{4.6}), we get
\begin{align}\label{4.8}
\|\y^n\phi_x(0,\cdot)\|_{H^1}&\lesssim\|\y^n[u^0_{xy}\psi+(v^0_p+v^1_e)u_{py}-u_{pyy}-F_p](0,\cdot)\|_{H^1}\nonumber\\
&\leq C_0\|\y^n\bar{u}_1\|_{H^3}+\|\y^nF_p(0,\cdot)\|_{H^1}.
\end{align}
In addition, the definition of $F_p$ gives
\[\|\y^nF_p(0,\cdot)\|_{H^3}\leq C_0,\]
substituting which into (\ref{4.8}) implies that
\begin{align}\label{4.9}
\|\y^n\phi_x(0,\cdot)\|_{H^1}\leq C_0\left(1+\|\y^n\bar{u}_1\|_{H^{3}}\right).
\end{align}

On the other hand, by the definition of $\phi$, we have
\[\psi=-u^0\int_y^\infty\frac{\phi(x,\theta)}{(u^0)^2}d\theta,\]
which implies that
\begin{align}\label{4.10}
v_p=-\psi_x=u^0_x\int_y^\infty\frac{\phi}{(u^0)^2}d\theta+u^0\int_y^\infty\left(\frac{\phi}{(u^0)^2}\right)_xd\theta.
\end{align}
Then, we get
\[\|\y^nv_p(0,\cdot)\|_{H^2}\leq C\left(\|\y^n\phi_y(0,\cdot)\|_{L^2}+\|\y^n\phi_x(0,\cdot)\|_{H^1}\right).\]
This estimate, together with (\ref{4.7}) and (\ref{4.9}), derives (\ref{4.4}).

Moreover, applying $\p_x$ to (\ref{4.6}) yields
\begin{align*}
\phi_{xx}=-u^0_{xxy}\psi+u^0_{xy}v_p-[v^0_{px}+v^1_{ex}]u_{py}+[v^0_p+v^1_e]v_{pyy}-v_{pyyy}+F_{px}.
\end{align*}
Then we obtain that
\begin{align}\label{4.11}
\|\y^n\phi_{xx}(0,\cdot)\|_{H^1}\leq C_0\|\y^n\bar{u}_1\|_{H^2}+C_0\|\y^nv_p(0,\cdot)\|_{H^4}+\|\y^nF_{px}(0,\cdot)\|_{H^1}.
\end{align}
It should be noted that
\begin{align}\label{4.12}
\|\y^nF_{px}(0,\cdot)\|_{H^1}\leq C_0(1+\|u^1_{exx}(0,\cdot)\|_{H^1}).
\end{align}
In addition, it follows from (\ref{4.2}) that
\begin{align}\label{4.13}
\|\y^nv_p(0,\cdot)\|_{H^4}&\leq \|\y^nG_p(0,\cdot)\|_{H^2}+\|\y^n\frac{u^0_{yy}v_p}{u^0}(0,\cdot)\|_{H^2}
+\|\y^n\frac{u_{py}}{u^0}(0,\cdot)\|_{H^4}\nonumber\\
&\leq C_0\left(\|\y^nF_p(0,\cdot)\|_{H^3}+\|\y^n\bar{u}_1\|_{H^5}+\|\y^nv_p(0,\cdot)\|_{H^2}\right)\nonumber\\
&\leq C_0(1+\|\y^n\bar{u}_1\|_{H^5}).
\end{align}
Therefore, in view of (\ref{4.10})-(\ref{4.13}), we have
\begin{align}\label{4.14}
\|\y^nv_{px}(0,\cdot)\|_{H^2}&\lesssim\|\y^n\phi_y(0,\cdot)\|_{L^2}+\|\y^n\phi_x(0,\cdot)\|_{H^1}
+\|\y^n\phi_{xx}(0,\cdot)\|_{H^1}\nonumber\\
&\leq C_0\left(1+\|\y^n\bar{u}_1\|_{H^5}+\|u^1_{exx}(0,\cdot)\|_{H^1}\right).
\end{align}
This completes the proof of this lemma.
\endProof

Step 2, we give the following auxiliary lemma.

\begin{lem}\label{L4.2}
For any $L>0$, denote $\Omega_\infty:=[0,L]\times \mathbb{R}_+$. Assume that $\p_x^jf,\p_x^jg\in L^2(\Omega_\infty),j=0,1$, decays fast as $y\rightarrow\infty$. Then the following fourth order partial deferential equation
\begin{align}\label{4.15}
-v_{xyy}+\frac{u^0_{yy}}{u^0}v_x+\left(\frac{v_{yy}}{u^0}\right)_{yy}=f_y+g
\end{align}
on $\Omega_\infty$ has an unique smooth solution satisfying initial data $v(0,y)=\bar{v}_0(y)$, boundary conditions $v=v_y=0$ on $y=0$ and $y=\infty$, and the estimate
\begin{align}\label{4.16}
&\sup_{0\leq x\leq L}\|\p_x^jv_{yy}(x)\|_{L^2}+\|\p_x^jv_{xy}\|_{L^2(\Omega_\infty)}\nonumber\\
&\leq C\sum^j_{i=0}\left(\|\p_x^iv_{yy}(0,\cdot)\|_{L^2}+\|\p_x^if\|_{L^2(\Omega_\infty)}+\|\y^\frac{3}{2}
\p_x^ig\|_{L^2(\Omega_\infty)}\right),j=0,1.
\end{align}
Moreover, there holds the weighted estimate
\begin{align}\label{4.161}
&\sup_{0\leq x\leq L}\|\y^n\p_x^jv_{yy}(x)\|_{L^2}+\|y^n\p_x^jv_{yyy}\|_{L^2(\Omega_\infty)}\nonumber\\
&\leq C\sum^j_{i=0}\left(\|\y^n\p_x^iv_{yy}(0,\cdot)\|_{L^2}+\|\y^n\p_x^if\|_{L^2(\Omega_\infty)}+\|\y^n
\p_x^ig\|_{L^2(\Omega_\infty)}\right),j=0,1.
\end{align}
\end{lem}

\Proof First, restrict the domain on $\Omega_N:=[0,L]\times[0,N]$ and give the approximate boundary conditions $v=v_y=0$ on $y=N$, instead of $y=\infty$. We introduce the inner product on $H^2(0,N)$:
\begin{align}\label{4.17}
[[u,v]]\equiv\int[u_yv_y+\frac{u^0_{yy}}{u^0}uv]dy
\end{align}
for any $u,v\in H^2(0,N)$. Let $\{e^i(y)\}_{i=1}^\infty$ be an orthogonal basis of $H^2(0,N)$ satisfying the same boundary conditions as $v$ doing. Here the orthogonality is obtained with respect to the inner product defined in (\ref{4.17}) and it holds that
\[[[e^i,e^j]]=\delta_{ij}, i,j\geq 1.\]
Base on (\ref{3.16}) and (\ref{3.17}), one can show that $[[\cdot,\cdot]]$ is equivalent to the usual inner product on $H^1(0,N)$. Then such an orthogonal basis exists.

 Now we introduce the weak formulation of (\ref{4.15}) as follows
\begin{align}\label{4.18}
[[v_x,e^i]]+\int\frac{v_{yy}e^i_{yy}}{u^0}dy=\int(-fv_{xy}+ge^i)dy
\end{align}
for any $e^i(y),i\geq 1.$ We will construct an approximate solution in Span$\{e^i(y)\}^k_{i=1}$ for (\ref{4.18}) defined as
\[v^k(x,y):=\sum_{j=1}^ka^j(x)e^j(y),\]
for each $k$. Substituting $v^k$ into (\ref{4.18}) in place of $v$, with the orthogonality of $\{e^i(y)\}_{i=1}^k$, yields
\begin{align}\label{4.19}
[[v^k_x,e^i]]+\int\frac{v^k_{yy}e^i_{yy}}{u^0}dy=\int(-fe^i_y+ge^i)dy,
\end{align}
 which is equivalent to a system of ODE equations:
\begin{align}\label{4.20}
a^i_x+\sum^k_{j=1}a^j\int\frac{e^j_{yy}e^i_{yy}}{u^0}dy=\int(-fe^i_y+ge^i)dy.
\end{align}
Since $f,g\in L^2(\Omega_N)$, there exists an unique smooth solution $(a^1,a^2,\cdots,a^k)$ for (\ref{4.20}), that is, there exists an unique smooth solution $v^k$ for (\ref{4.18}). In order to take $k$ tends to infinity, we need some energy estimates.

Multiplying (\ref{4.19}) by $a_x^i$ and take the sum over $i$ from $1$ to $k$, we get
\begin{align}\label{4.21}
[[v^k_x,v^k_x]]+\frac{1}{2}\frac{d}{dx}\int\frac{(v^k_{yy})^2}{u^0}dy=\int\left(-fv^k_{xy}+gv^k_{x}
+\frac{1}{2}\left(\frac{1}{u^0}\right)_x(v^k_{yy})^2\right)dy.
\end{align}
Similar to the analysis in (\ref{3.14}) and (\ref{3.15}), we have
\begin{align}\label{4.22}
[[v^k_{x},v^k_{x}]]\geq \alpha\|v^k_{xy}\|^2_{L^2(0,N)}.
\end{align}
Then, applying the Gronwall's inequality gives that
\begin{align}\label{4.23}
\sup_{x\in[0,L]}\|v^k_{yy}\|^2_{L^2(0,N)}+\|v^k_{xy}\|^2_{L^2(\Omega_N)}
\lesssim\|v^k_{yy}(0,\cdot)\|^2_{L^2(0,N)}+\|(f,\y^{\frac{3}{2}}g)\|_{L^2(\Omega_N)}^2.
\end{align}
Taking $k\rightarrow\infty$ yields the weak solution $v(x,y)$ to (\ref{4.15}), which satiesfies
\begin{align}\label{4.231}
\sup_{x\in[0,L]}\|v_{yy}\|^2_{L^2(0,N)}+\|v_{xy}\|^2_{L^2(\Omega_N)}
\lesssim\|v_{yy}(0,\cdot)\|^2_{L^2(0,N)}+\|(f,\y^\frac{3}{2}g)\|_{L^2(\Omega_N)}^2.
\end{align}

Next, we should derive higher regularity for the weak solution. Applying $\p_x$ to (\ref{4.19}), multiplying the result by $a^i_{xx}$ and take the sum over $i$ from $1$ to $k$, we have
\begin{align}\label{4.24}
&[[v^k_{xx},v^k_{xx}]]+\frac{1}{2}\frac{d}{dx}\int\frac{(v^k_{xyy})^2}{u^0}dy\nonumber\\
=&\frac{1}{2}\int\left(\frac{1}{u^0}\right)_x(v^k_{xyy})^2dy-\int\left(\frac{u^0_{yy}}{u^0}\right)_xv^k_xv^k_{xx}dy
-\int\left(\frac{1}{u^0}\right)_xv^k_{yy}v^k_{xxyy}dy\nonumber\\
&+\int\left(-f_xv^k_{xxy}+g_xv^k_{xx}\right)dy:=\mathcal{J}_1+\mathcal{J}_2+\mathcal{J}_3+\mathcal{J}_4.
\end{align}
Note that
\begin{align}\label{4.25}
\mathcal{J}_1+\mathcal{J}_3&=-\frac{d}{dx}\int\left(\frac{1}{u^0}\right)_xv^k_{yy}v^k_{xyy}dy
+\int\left(\frac{1}{u^0}\right)_{xx}v^k_{yy}v^k_{xyy}dy+3\mathcal{J}_1\nonumber\\
&\leq -\frac{d}{dx}\int\left(\frac{1}{u^0}\right)_xv^k_{yy}v^k_{xyy}+C(u^0)(\|v^k_{xyy}\|^2_{L^2(0,N)}+\|v^k_{yy}\|^2_{L^2(0,N)}),
\end{align}
and that
\begin{align}\label{4.26}
\mathcal{J}_2+\mathcal{J}_4&\leq \delta\|v^k_{xxy}\|^2_{L^2(0,N)}+C(u^0)\left(\|v^k_{xy}\|^2_{L^2(0,N)}+\|(f_x,\y^\frac{3}{2}g_x)\|^2_{L^2(0,N)}\right).
\end{align}
In addition, similar to (\ref{4.22}), we have
\begin{align}\label{4.27}
[[v^k_{xx},v^k_{xx}]]\geq \alpha\|v^k_{xxy}\|^2_{L^2(\Omega_N)}.
\end{align}
Then substituting (\ref{4.25})-(\ref{4.27}) into (\ref{4.26}) and applying the Gronwall's inequality gives
\begin{align*}
\sup_{0\leq x\leq L}\|v^k_{xyy}\|^2_{L^2(0,N)}+\|v^k_{xxy}\|^2_{L^2(\Omega_N)}\lesssim& \sum^1_{i=0}\|\p^i_xv^k_{yy}(0,\cdot)\|^2_{L^2(0,N)}+\sup_{0\leq x\leq L}\|v^k_{yy}\|^2_{L^2(0,N)}\\
&+\|v^k_{xy}\|^2_{L^2(\Omega_N)}+\|(f_x,\y^\frac{3}{2}g_x)\|^2_{L^2(\Omega_N)},
\end{align*}
which together with (\ref{4.23}) gives
 \begin{align}\label{4.28}
&\sup_{x\in[0,L]}\|v^k_{xyy}\|^2_{L^2(0,N)}+\|v^k_{xxy}\|^2_{L^2(\Omega_N)}\nonumber\\
&\lesssim\sum^1_{i=0}\left(\|\p_x^iv^k_{yy}(0,\cdot)\|^2_{L^2(0,N)}+\|(\p_x^if,\y^\frac{3}{2}\p_x^ig)\|_{L^2(\Omega_N)}^2\right)
 \end{align}
Again, taking $k\rightarrow\infty$ yields
\begin{align}\label{4.281}
&\sup_{x\in[0,L]}\|\p^j_xv_{yy}\|^2_{L^2(0,N)}+\|\p_x^jv_{xy}\|^2_{L^2(\Omega_N)}\nonumber\\
&\lesssim\sum^j_{i=0}\left(\|\p_x^iv_{yy}(0,\cdot)\|^2_{L^2(0,N)}+\|(\p_x^if,\y^\frac{3}{2}\p_x^ig)\|_{L^2(\Omega_N)}^2\right),j=0,1.
 \end{align}
It should be noted that all the constants $C$ in the estimates above are independent of $N$, and hence the unweighted estimates (\ref{4.16}) is proved as taking $N\rightarrow \infty$.

Finally, let us derive the weighted estimates. The readers should notice that the weight function for diffusion terms is $y^n$, but we sometimes write the other terms by weight function $\y^n$ since $y^n\leq \y^n$.

On one hand, multiplying (\ref{4.15}) by $y^{2n}v_{yy}$ and integrating by part over $\mathbb{R}_+$, we get
\begin{align}\label{4.282}
&\frac{1}{2}\frac{d}{dx}\int |y^nv_{yy}|^2+\int\frac{|y^nv_{yyy}|^2}{u^0}\nonumber\\
=&\int\frac{u^0_{pyy}}{u^0}v_xv_{yy}y^{2n}+\int\left(\frac{v_{yy}}{u^0}\right)_yv_{yy}(y^{2n})_y
+\int\frac{u^0_{py}}{(u^0)^2}v_{yy}v_{yyy}y^{2n}\nonumber\\
&+\int fv_{yyy}y^{2n}+\int fv_{yy}(y^{2n})_y+\int gv_{yy}y^{2n}:=\sum_{i=1}^6\mathcal{K}_i.
\end{align}
Note that
\begin{align*}
\mathcal{K}_1\lesssim&\int|u^0_{pyy}||v_x||v_{yy}|y^{2n}\lesssim\|v_{xy}\|_{L^2}\int|u^0_{pyy}||v_{yy}|\y^{2n+1}\\
\lesssim&\|v_{xy}\|_2\|\y^{n+1}u^0_{pyy}\|_2\|\y^n v_{yy}\|_2\leq C\|\y^n v_{yy}\|^2_2+C\|v_{xy}\|^2_2,\\
\mathcal{K}_2+\mathcal{K}_3\lesssim&\int|v_{yyy}||v_{yy}|y^{2n}+\int|v_{yy}|^2y^{2n}\leq \delta\|y^n v_{yyy}\|^2_2
+C\|\y^n v_{yy}\|^2_2,\\
\sum_{i=4}^6\mathcal{K}_i\lesssim&\int|f||v_{yyy}|y^{2n}+\int|f||v_{yy}|y^{2n}+\int|g||v_{yy}|y^{2n}\\
\leq& \delta\|y^nv_{yyy}\|^2_2+C\|\y^nv_{yy}\|^2_2+C\|\y^nf\|^2_2+C\|\y^ng\|^2_2.
\end{align*}
Substituting these estimates into (\ref{4.282}) with taking $\delta$ small enough yields
\begin{align*}
&\frac{d}{dx}\|y^nv_{yy}\|^2_2+\|y^nv_{yyy}\|^2_2\lesssim \|\y^nv_{yy}\|^2_2+\|v_{xy}\|^2_2+\|\y^nf\|^2_2+\|\y^ng\|^2_2,
\end{align*}
applying Gronwall's inequality to which, together with estimate (\ref{4.16}), implies that
\begin{align}\label{4.283}
\sup_{x\in[0,L]}\|\y^nv_{yy}\|^2_2+\|y^nv_{yyy}\|^2_{L^2(\Omega_\infty)}
\lesssim\|\y^nv_{yy}(0,\cdot)\|^2_2+\|\y^n(f,g)\|^2_{L^2(\Omega_\infty)}.
\end{align}

On the other hand, applying $\p_x$ to equation (\ref{4.15}), multiplying the result by $y^{2n}v_{xyy}$ and integrating by part over $\mathbb{R}_+$, we have
\begin{align}\label{4.284}
&\frac{1}{2}\frac{d}{dx}\int|y^nv_{xyy}|^2+\int\frac{|y^nv_{xyyy}|^2}{u^0}\nonumber\\
=&\int\left(\frac{u^0_{pyy}v_x}{u^0}\right)_xv_{xyy}y^{2n}-\int\frac{v_{yyy}u^0_{px}}{(u^0)^2}v_{xyyy}y^{2n}
-\int\frac{v_{xyy}u^0_{py}}{(u^0)^2}v_{xyyy}y^{2n}\nonumber\\
&-\int\left(\frac{1}{u^0}\right)_{xy}v_{yy}v_{xyyy}y^{2n}-\int\left(\frac{v_{yy}}{u^0}\right)_{xy}v_{xyy}(y^{2n})_y\nonumber\\
&+\int f_x(v_{xyy}y^{2n})_y-\int g_xv_{xyy}y^{2n}:=\sum_{i=1}^7\mathcal{L}_i.
\end{align}

Similar to the estimates on $\mathcal{K}_i$, we infer that
\begin{align*}
\mathcal{L}_1&\lesssim\int(|u^0_{pxyy}||v_x|+|u^0_{pyy}||u^0_{px}||v_x|+|u^0_{pyy}||v_{xx}|)|v_{xyy}|y^{2n}\\
&\lesssim \|v_{xy}\|^2_2+\|v_{xxy}\|^2_2+\|\y^nv_{xyy}\|^2_2,\\
\sum_{i=2}^4\mathcal{L}_i&\lesssim\int(|v_{yyy}||u^0_{px}|+|v_{xyy}||u^0_{py}|+|v_{yy}||u^0_{pxy}|+|v_{yy}||u^0_{px}||u^0_{py}|)
|v_{xyyy}|y^{2n}\\
&\leq \delta\|y^nv_{xyyy}\|^2_2+C\|y^nv_{yyy}\|^2_2+C\left(\|\y^nv_{xyy}\|^2_2+\|\y^nv_{yy}\|^2_2\right),\\
\mathcal{L}_5&\lesssim\int(|u^0_x||u^0_y||v_{yy}|+|u^0_{xy}||v_{yy}|+|u^0_x||v_{yyy}|+|v_{xyyy}|+|u^0_y||v_{xyy}|)|v_{xyy}|y^{2n-1}\\
&\leq \delta\|y^nv_{xyyy}\|^2_2+C\|y^nv_{yyy}\|^2_2+C\left(\|\y^nv_{xyy}\|^2_2+\|\y^nv_{yy}\|^2_2\right),\\
\sum_{i=6}^7\mathcal{L}_i&\lesssim\int|f_x||v_{xyyy}|y^{2n}+\int|f_x||v_{xyy}|y^{2n}+\int|g_x||v_{xyy}|y^{2n}\\
&\leq \delta\|y^nv_{xyyy}\|^2_2+C\left(\|\y^nv_{xyy}\|^2_2+\|\y^n(f_x,g_x)\|^2_2\right).
\end{align*}
Substituting the estimates of $\mathcal{L}_i$ into (\ref{4.284}) with $\delta$ sufficiently small, we get
\begin{align}\label{4.285}
&\frac{d}{dx}\|y^nv_{xyy}\|_2^2+\|y^nv_{xyyy}\|^2_2\nonumber\\
&\lesssim \sum_{i=0}^1(\|\y^n\p_x^iv_{yy}\|_2^2+\|\p_x^iv_{xy}\|_2^2)+\|y^nv_{yyy}\|_2^2+\|\y^n(,f_x,g_x)\|^2_2.
\end{align}
Finally, applying Gronwall's inequality with using estimates (\ref{4.16}) and (\ref{4.283}) derive (\ref{4.161}).
\endProof

Step 3, with these two lemmas in hand, we are able to prove the existence of smooth solutions for system (\ref{4.1}) on $\Omega_\infty:=[0,L]\times \mathbb{R}_+$.

\begin{lem}\label{L4.3}
Under the assumptions in Theorem 1.1, there exists an unique smooth solution $[u_p,v_p]$ to system (\ref{4.1}) satisfying that
\begin{align}\label{4.29}
\|[u_p,v_p]\|_{L^\infty(\Omega_\infty)}+\sup_{0\leq x\leq L}\|\y^nv_{pyy}\|_{L^2(\mathbb{R}_+)}+\|y^nv_{pyyy}\|_{L^2(\Omega_\infty)}\leq C(L,\kappa)\e^{-\kappa},
\end{align}
for some $\kappa>0$ small sufficiently. Moreover, the following higher regularity estimate
\begin{align}\label{4.30}
\sup_{0\leq x\leq L}\|\y^nv_{pxyy}\|_{L^2(\mathbb{R}_+)}+\|y^nv_{pxyyy}\|_{L^2(\Omega_\infty)}\leq C(L)\e^{-1}
\end{align}
holds uniformly in small $\e$, in which the constant $C(L)$ depends only on $[u^0,v^0_p]$, the given boundary data and $L$.
\end{lem}

\Proof Denote
\begin{align}\label{4.31}
\bar{v}&:=v_p-y\chi(y)u^1_{ex}(x,0):=v_p+\tilde{v}.
\end{align}
Then $\bar{v}=\bar{v}_y=0$ on $y=0$ and $y=\infty.$ Furthermore, $\bar{v}$ satisfies equation
\begin{align}\label{4.32}
-\bar{v}_{xyy}+\frac{u^0_{yy}}{u^0}\bar{v}_x+\left\{\frac{\bar{v}_{yy}}{u^0}\right\}_{yy}
=&G_{px}-\left\{\frac{u^0_{yy}}{u^0}\right\}_xv_p+\left(u_{py}\left\{\frac{1}{u^0}\right\}_x\right)_{yy}\nonumber\\
&+\tilde{v}_{xyy}-\frac{u^0_{yy}}{u^0}\tilde{v}_x-\left\{\frac{\tilde{v}_{yy}}{u^0}\right\}_{yy}:=f_y+g,
\end{align}
where we denote
\begin{align}
f:=&\left(u_{py}\left\{\frac{1}{u^0}\right\}_x\right)_y+\frac{v^0_p+v^1_e}{u^0}v_{pyy}
+2v_{pyy}\left\{\frac{1}{u^0}\right\}_y-\frac{v_{ex}^1 u_{py}}{u^0}+\frac{F_{px}}{u^0}=\sum_{i=1}^5f_i,\label{4.33}\\
g:=&\left\{\frac{1}{u^0}\right\}_x\left[F_{py}-u_pu^0_{xy}+(v^0_p+v^1_e)u_{pyy}-\sqrt\e v^1_{ez}u_{py}\right]+\left\{\frac{v^1_{ex}}{u^0}\right\}_yu_{py}\nonumber\\
&-\left(\frac{v^0_p+v^1_e}{u^0}\right)_yv_{pyy}-2u_{pyy}\left\{\frac{1}{u^0}\right\}_{xy}-u_{py}\left\{\frac{1}{u^0}\right\}_{xyy}
-v_{pyy}\left\{\frac{1}{u^0}\right\}_{yy}\nonumber\\
&-F_{px}\left(\frac{1}{u^0}\right)_y-\frac{1}{u^0}(u_{px}u^0_{xy}+u_pu^0_{xxy}+v^0_{px}u_{pyy})-\frac{\sqrt\e}{u^0}(v^1_{exz}u_{py}+ v^1_{ez}u_{pxy})\nonumber\\
&-\left\{\frac{u^0_{yy}}{u^0}\right\}_xv_p+\tilde{v}_{xyy}-\frac{u^0_{yy}}{u^0}\tilde{v}_x
-\left\{\frac{\tilde{v}_{yy}}{u^0}\right\}_{yy}=\sum_{i=1}^{13}g_i,\label{4.34}
\end{align}
with
\[F_p=-u^0_{px}u^1_e-u^0_pu^1_{ex}-(v^0_p+yu^0_{px})u^0_{ez}-yu^0_{py}v^1_{ez}+E_1.\]
Due to the divergence-free condition, we infer that $u_{px}=-v_{py}=-\bar{v}_y+\tilde{v}_y$. We shall work with the norm:
\begin{align}\label{4.35}
|||\bar{v}|||^2\equiv\sup_{0\leq x\leq L}\|\y^n\bar{v}_{yy}(x)\|^2_2+\|y^n\bar{v}_{yyy}\|^2_{L^2(\Omega_\infty)}.
\end{align}
In view of Lemma \ref{L4.2}, we have
\begin{align}\label{4.36}
|||\bar{v}|||^2\leq C\left(\|\y^n\bar{v}_{yy}(0,\cdot)\|^2_2+\|\y^nf\|^2_{L^2(\Omega_\infty)}+\|\y^ng\|^2_{L^2(\Omega_\infty)}\right)
\end{align}
with $(f,g)$ being defined as in (\ref{4.33}) and (\ref{4.34}). Recall that $\bar{v}=v_p+\tilde{v}$. Then  by the definition of $\tilde{v}$ and estimate (\ref{4.4}), it follows that
\begin{align}\label{4.37}
\|\y^n\bar{v}_{yy}(0,\cdot)\|_2&\leq \|\y^nv_{pyy}(0,\cdot)\|_2+C\leq C_0(1+\|\y^n\bar{u}_1\|_{H^3}),
\end{align}
where we have used the fact that
\[|u^1_{ex}(0,0)|\leq \int_0^1|[(z-1)v^1_{ez}(0,z)]_z|dz\leq \|V_{b0}\|_{H^2(0,1)}\leq C.\]
Next, let us give bounds on $f$ and $g$. We infer that
\begin{align}
&|\bar{v}_y|\leq \int_y^\infty|\bar{v}_{yy}|d\theta\lesssim\|\y^n\bar{v}_{yy}\|_2\y^{-n+1},\label{4.38}\\
&|\bar{v}|\leq\int_y^\infty|\bar{v}_y|d\theta\lesssim \|\y^n\bar{v}_{yy}\|_2\y^{-n+3}.\label{4.39}
\end{align}
Thus, we obtain
\begin{align*}
\int_0^\infty|\bar{v}|^2&\leq C\|\y^n\bar{v}_{yy}\|^2_2\int_0^\infty\y^{-2n+6}\leq C|||\bar{v}|||^2,\\
\int_0^\infty|\bar{v}_y|^2&\leq C\|\y^n\bar{v}_{yy}\|^2_2\int_0^\infty\y^{-2n+2}\leq C|||\bar{v}|||^2,
\end{align*}
for some $n$ large enough. In addition, there hold that
\begin{align*}
\int_0^\infty|u_p|^2\leq& \int_0^\infty|\bar{u}_1|^2+L\int_0^\infty\int_0^L|u_{px}|^2\\
\leq& \|\bar{u}_1\|^2_2+L\int_{\Omega_\infty}|\bar{v}_y|^2+L\int_{\Omega_\infty}|\tilde{v}_y|^2\\
\leq &CL|||\bar{v}|||^2+C\left(\|\y^n\bar{u}_1\|^2_2+\|u^1_{ex}(x,0)\|^2_{L^2(0,L)}\right),\\
\int_0^\infty\y^{2n}|u_{py}|^2\leq& \|\y^n\bar{u}^\prime_1\|^2_2+L\int_{\Omega_\infty}\y^{2n}|u_{pxy}|^2\\
\leq& CL|||\bar{v}|||^2+C\left(\|\y^n\bar{u}^\prime_1\|^2_2+\|u^1_{ex}(x,0)\|^2_{L^2(0,L)}\right),
\end{align*}
Hence, in view of equation (\ref{4.1}), we get
\begin{align*}
&\iint_{\Omega_\infty}|u_{pyy}|^2\leq \int_{\Omega_\infty}|u^0u_{px}+u_pu^0_x+v_pu^0_y+[v^0_p+v^1_e]u_{py}-F_p|^2\\
\leq &\iint_{\Omega_\infty}\left(|\bar{v}_y|^2+|\tilde{v}_y|^2+|u_p|^2+|\bar{v}|^2+|\tilde{v}|^2
+|u_{py}|^2+|F_p|^2\right)\\
\leq &CL|||\bar{v}|||^2+C\left(\|\y^n\bar{u}_1\|^2_{H^1}+\|u^1_{ex}(x,0)\|^2_{L^2(0,L)}+1\right)
\end{align*}
Therefore, we have
\begin{align}
\iint_{\Omega_\infty}\y^{2n}|f_1+f_4|^2&\lesssim \iint_{\Omega_\infty}|u_{pyy}|^2+\int_0^L\|v^1_{ex}\|^2_\infty\int_0^\infty\y^{2n}|u_{py}|^2\nonumber\\
&\leq CL|||\bar{v}|||^2+C\left(\|\y^n\bar{u}_1\|^2_{H^1}+\|u^1_{ex}(x,0)\|^2_{L^2(0,L)}+1\right)\label{4.40}\\
\iint_{\Omega_\infty}\y^{2n}|f_2+f_3|^2
&\lesssim\iint_{\Omega_\infty}\y^{2n}|v_{pyy}|^2\leq CL|||\bar{v}|||^2+C\|u^1_{ex}(x,0)\|^2_{L^2(0,L)}.\label{4.41}\\
\iint_{\Omega_\infty}\y^{2n}|f_5|^2&\lesssim\int_0^L(\|u^1_e\|^2_\infty+\|v^1_{ez}\|^2_\infty+\|v^1_{exz}\|^2_\infty)+1
\nonumber\\
&\leq C(1+\|v^1_e\|_{W^{3,q}})\leq C(L,\kappa)\e^{-2\kappa}\label{4.411}
\end{align}
Similarly, it follows that
\begin{align}
\iint_{\Omega_\infty}\y^{2n}|g_1|^2&\lesssim \iint_{\Omega_\infty}(|F_{py}|^2+|u_p|^2+|u_{pyy}|^2)+\e\int_0^L\|v^1_{ez}\|^2_\infty\int_0^\infty|u_{py}|^2\nonumber\\
&\leq CL|||\bar{v}|||^2+C\left(\|\y^n\bar{u}_1\|^2_{H^1}+\|u^1_{ex}(x,0)\|^2_{L^2(0,L)}+1\right),\label{4.42}\\
\iint_{\Omega_\infty}\y^{2n}|g_2|^2&\lesssim\int_0^L(\e\|v^1_{exz}\|^2_\infty+\|v^1_{ex}\|^2_\infty)
\int_0^\infty\y^{2n}|u_{py}|^2\nonumber\\
&\leq CL|||\bar{v}|||^2+C\left(\|\y^n\bar{u}_1\|^2_{H^1}+\|u^1_{ex}(x,0)\|^2_{L^2(0,L)}+1\right),\label{4.421}\\
\iint_{\Omega_\infty}\y^{2n}|g_3|^2&\lesssim\left(\e\|v^1_{ez}\|^2_\infty+\|v^1_e\|_\infty+1\right)
\iint_{\Omega_\infty}\y^{2n}|v_{pyy}|^2\nonumber\\
&\leq CL|||\bar{v}|||^2+C\left(\|u^1_{ex}(x,0)\|^2_{L^2(0,L)}+1\right),\label{4.422}\\
\sum^6_{i=4}\iint_{\Omega_\infty}\y^{2n}|g_i|^2&\lesssim\iint_{\Omega_\infty}\y^{2n}(|u_{py}|^2+|v_{pyy}|^2)
+\iint_{\Omega_\infty}|u_{pyy}|^2\nonumber\\
&\leq CL|||\bar{v}|||^2+C\left(\|\y^n\bar{u}_1\|^2_{H^1}+\|u^1_{ex}(x,0)\|^2_{L^2(0,L)}+1\right),\label{4.43}\\
\iint_{\Omega_\infty}\y^{2n}|g_7|^2&\lesssim\iint_{\Omega_\infty}\y^{2n}|f_5|^2\leq C(L,\kappa)\e^{-2\kappa},\label{4.431}\\
\iint_{\Omega_\infty}\y^{2n}|g_8|^2&\lesssim\iint_{\Omega_\infty}(|v_{py}|^2+|u_p|^2+|u_{pyy}|^2)\nonumber\\
&\leq CL|||\bar{v}|||^2+C\left(\|\y^n\bar{u}_1\|^2_{H^1}+\|u^1_{ex}(x,0)\|^2_{L^2(0,L)}+1\right),\label{4.44}\\
\iint_{\Omega_\infty}\y^{2n}|g_9|^2&\lesssim\e\int_0^L\|v^1_{exz}\|^2_{L^\infty}\int_0^\infty\y^{2n}|u_{py}|^2
+\e\|v^1_{ez}\|^2_{L^\infty(\Omega_0)}\iint_{\Omega_\infty}\y^{2n}|v_{pyy}|^2\nonumber\\
&\leq CL|||\bar{v}|||^2+C\left(\|\y^n\bar{u}_1\|^2_{H^1}+\|u^1_{ex}(x,0)\|^2_{L^2(0,L)}+1\right),\label{4.441}\\
\sum^{13}_{i=10}\iint_{\Omega_\infty}\y^{2n}|g_i|^2&\lesssim\iint_{\Omega_\infty}|\bar{v}|^2+\int^L_0(|u^1_{ex}(x,0)|^2
+|u^1_{exx}(x,0)|^2)\nonumber\\
&\leq CL|||\bar{v}|||^2+C\|u^1_{ex}(x,0)\|^2_{L^2(0,L)}+C\|u^1_{exx}(x,0)\|^2_{L^2(0,L)}.
\end{align}

It should be noted that
\begin{align}\label{4.45}
&\sum_{i=0}^1\|\p_x^iu^1_{ex}(x,0)\|^2_{L^2(0,L)}=\sum_{i=0}^1\|\p_x^iv^1_{ez}(x,0)\|^2_{L^2(0,L)}\nonumber\\
&\leq \sum_{i=0}^1\iint_{\Omega_0}|[(z-1)(\p_x^iv^1_{ez})^2]_z|
\leq \sum_{i=0}^1\iint_{\Omega_0}|\p_x^iv^1_{ez}|^2+\sum_{i=0}^1\iint_{\Omega_0}|\p_x^iv^1_{ez}||\p_x^iv^1_{ezz}|\nonumber\\
&\leq \sum_{i=0}^1\|\p_x^iv^1_{ez}\|^2_2+\sum_{i=0}^1\|\p_x^iv^1_{ez}\|_p\|\p_x^iv^1_{ezz}\|_q
\leq C\e^{-2\kappa},
\end{align}
for sufficiently small $\kappa>0$. In conclusion, we obtain
\begin{align}\label{4.46}
|||\bar{v}|||\leq CL|||\bar{v}|||+C(L,\kappa)\e^{-\kappa},
\end{align}
which with sufficiently small $L$ give the uniform bound for $|||\bar{v}|||$.

Furthermore, since equation (\ref{4.32}) is linear with respect to $\bar{v}$, together with estimate (\ref{4.46}), it is easy to apply the contraction mapping theorem to show the existence of the unique solution for (\ref{4.32}) and hence (\ref{4.3}). Then, it follows from the boundedness of $|||\bar{v}|||$ that
\begin{align}\label{4.461}
\sup_{0\leq x\leq L}\|\y^nv_{pyy}\|^2_2+\int_0^L\|y^nv_{pyyy}\|^2_2\leq C(L,\kappa)\e^{-2\kappa}.
\end{align}

The boundedness of $v_p$ follows by the calculation similar to (\ref{4.38}) and (\ref{4.39}):
\begin{align}\label{4.47}
|v_p(x,y)|&\leq \int_y^\infty|v_{py}|d\theta\leq \int_{\mathbb{R}_+}\y^{-n+1}\|\y^nv_{pyy}\|_2dy\leq C\|\y^nv_{pyy}\|_2,
\end{align}
which implies that
\[\|v_p\|_\infty\lesssim\sup_{0\leq x\leq L}\|\y^nv_{pyy}\|_2\leq C(L,\kappa)\e^{-\kappa}.\]
Similarly, the boundedness of $u_p$ follows from the definition
\[u_p(x,y)=\bar{u}_1(y)-\int_0^xv_{py}ds,\]
which gives that
\[|u_p(x,y)|\leq |\bar{u}_1(y)|+\int_0^L|v_{py}|dx\leq \|\y^n\bar{u}_1^\prime\|_2+C\int_0^L\|\y^nv_{pyy}\|_2.\]
Then, we have
\[\|u_p\|_\infty\lesssim \|\y^n\bar{u}_1^\prime\|_2+\sup_{0\leq x\leq L}\|\y^nv_{pyy}\|_2\leq C(L,\kappa)\e^{-\kappa}.\]

To complete the proof of the lemma, we are now concerned with the higher regularity estimate. Again, applying Lemma \ref{L4.2} to
equation (\ref{4.32}), we also get
\begin{align*}
&\sup_{0\leq x\leq L}\|\y^n\bar{v}_{xyy}\|^2_2+\int_0^L\|y^n\bar{v}_{xyyy}\|^2_2\nonumber\\
&\leq C\sum^1_{i=0}\left(\|\y^n\p_x^i\bar{v}_{yy}(0,\cdot)\|^2_2+\int_0^L\|\y^n\p_x^if\|^2_2+\int_0^L\|\y^n
\p_x^ig\|^2_2\right),
\end{align*}
which, in view of estimate (\ref{4.46}), is reduced to
\begin{align}\label{4.48}
&\sup_{0\leq x\leq L}\|\y^n\bar{v}_{xyy}\|^2_2+\int_0^L\|y^n\bar{v}_{xyyy}\|^2_2\nonumber\\
&\lesssim\|\y^n\bar{v}_{xyy}(0,\cdot)\|^2_2+\int_0^L\|\y^nf_x\|^2_2+\int_0^L\|\y^ng_x\|^2_2+C(L,\kappa)\e^{-2\kappa}.
\end{align}
Recalled by (\ref{4.5}) that
\begin{align*}
\|\y^n\bar{v}_{xyy}(0,\cdot)\|_2&\leq \|\y^nv_{pxyy}(0,\cdot)\|_2+C|u^1_{exx}(0,0)|\\
&\lesssim 1+|u^1_{exx}(0,0)|+\|u^1_{exx}(0,\sqrt\e\cdot)\|_{H^1}.
\end{align*}
Note that
$$Lf^2(0)=-\int_0^L\p_x[(L-x)f^2(x)]dx\leq \|f\|^2_{L^2}+2L\|f\|_{L^2}\|f_x\|_{L^2},$$
which gives
$$|f(0)|\leq L^{-1/2}\|f\|_{L^2(0,L)}+\sqrt 2\|f\|^{1/2}_{L^2(0,L)}\|f_x\|^{1/2}_{L^2(0,L)}.$$
Then, using the estimate of $v^1_e$ in Section 3, we have
\begin{align*}
\|u^1_{exx}(0,\sqrt\e\cdot)\|_{L^2(\mathbb{R}_+)}&\leq C\e^{-1/4}\|v^1_{exz}(0,\cdot)\|_{L^2(0,1)}\\
&\leq C(L)\e^{-1/4}\left(\|v^1_{exz}\|_{L^2(\Omega_0)}+\|v^1_{ezx}\|^{1/2}_{L^2(\Omega_0)}
\|v^1_{exxz}\|^{1/2}_{L^2(\Omega_0)}\right)\\
&\leq C(L)\e^{-1/4}\e^{-1/4}\leq C(L)\e^{-1/2},\\
\|u^1_{exxy}(0,\sqrt\e\cdot)\|_{L^2(\mathbb{R}_+)}&\leq C\e^{1/4}\|v^1_{exzz}(0,\cdot)\|_{L^2(0,1)}\\
&\leq C(L)\e^{1/4}\left(\|v^1_{exzz}\|_{L^2(\Omega_0)}+\|v^1_{exzz}\|^{1/2}_{L^2(\Omega_0)}
\|v^1_{exxzz}\|^{1/2}_{L^2(\Omega_0)}\right)\\
&\leq C(L)\e^{1/4}\e^{-1}\leq C(L)\e^{-3/4}.
\end{align*}
Also, there holds
\begin{align*}
|u^1_{exx}(0,0)|^2&\leq \|u^1_{exx}(0,\cdot)\|_{L^2(0,1)}\|u^1_{exxz}(0,\cdot)\|_{L^2(0,1)}\\
&\leq \|v^1_{exz}(0,\cdot)\|_{L^2(0,1)}\|v^1_{exzz}(0,\cdot)\|_{L^2(0,1)}\\
&\leq C(L)\e^{-1/4}\e^{-1}\leq C(L)\e^{-3/2}.
\end{align*}
These implies that
\[\|\y^n\bar{v}_{xyy}(0,\cdot)\|_2\leq C(L)\e^{-3/4},\]
and hence it follows from (\ref{4.48}) that
\begin{align}\label{4.49}
&\sup_{0\leq x\leq L}\|\y^nv_{pxyy}\|^2_2+\int_0^L\|y^nv_{pxyyy}\|^2_2\nonumber\\
&\lesssim C(L)\e^{-3/2}+\sup_{0\leq x\leq L}\|u^1_{exx}(x,0)\|^2+\|u^1_{exx}(x,0)\|^2_{L^2(0,L)}+\int_0^L\|\y^n(f_x,g_x)\|^2_2.
\end{align}
We infer that
\begin{align*}
\|u^1_{exx}(x,0)\|^2_{L^2(0,L)}&\leq 2\int_{\Omega_0}|u^1_{exx}u^1_{exxz}|\leq 2\int_{\Omega_0}|v^1_{exz}v^1_{exzz}|\\
&\leq C\|v^1_{exz}\|_{L^2(\Omega_0)}\|v^1_{exzz}\|_{L^2(\Omega_0)}\leq C\e^{-1/2},\\
\|u^1_{exxx}(x,0)\|^2_{L^2(0,L)}&\leq 2\int_{\Omega_0}|u^1_{exxx}u^1_{exxxz}|\leq 2\int_{\Omega_0}|v^1_{exzz}v^1_{exxzz}|\\
&\leq C\|v^1_{exxz}\|_{L^2(\Omega_0)}\|v^1_{exxzz}\|_{L^2(\Omega_0)}\leq C\e^{-2},
\end{align*}
and then
\begin{align*}
\sup_{0\leq x\leq L}\|u^1_{exx}(x,0)\|_2^2\leq L\|u^1_{exxx}(x,0)\|^2_{L^2(0,L)}+|u^1_{exx}(0,0)|^2\leq C\e^{-2}.
\end{align*}
Substituting these estimates for boundary terms into (\ref{4.49}), we obtain
\begin{align}\label{4.50}
\sup_{0\leq x\leq L}\|\y^nv_{pxyy}\|^2_2+\int_0^L\|y^nv_{pxyyy}\|^2_2
\leq C(L)\e^{-2}+C\int_0^L\|\y^n(f_x,g_x)\|^2_2.
\end{align}
Since that the estimates for $f_x,g_x$ are similarly as done above, we omit the details here. The proof of this lemma is completed.
\endProof

Since we will use the estimates on $v^1_{px},v^1_{pxy}$, and $v^1_{pxx}$ in estimating the $L^2$ norm of $R^u_{app}$ and $R^v_{app}$, we give the following
\begin{cor}
Let $v_{p}$ be the solution constructed in Lemma \ref{L4.3}, then it follows that
\begin{align}\label{4.50a}
\|\y^n(v_{px},v_{pxy})\|_{L^2(\Omega_\infty)}\leq C(L,\kappa)\e^{-\kappa},\ \ \|\y^n(v_{pxx},v_{pxxy})\|_{L^2(\Omega_\infty)}\leq C(L)\e^{-1},
\end{align}
for any $n\in\mathbb{N}_+.$
\end{cor}
\Proof
By virtue of (\ref{4.16}), we have
\begin{align*}
\|\p_x^j\bar{v}_{xy}\|_{L^2(\Omega_\infty)}&\lesssim \sum^j_{i=0}\left(\|\p_x^i\bar{v}_{yy}(0,\cdot)\|_2+\|(\p^j_xf,\y^{3/2}\p^j_xg)\|_{L^2(\Omega_\infty)}\right),\\
&\lesssim \sum^j_{i=0}\left(\|\y^n\p_x^i\bar{v}_{yy}(0,\cdot)\|_2+\|\y^n(\p^j_xf,\p^j_xg)\|_{L^2(\Omega_\infty)}\right),~~j=0,1,
\end{align*}
where $f$ and $g$ are defined as in (\ref{4.33}) and (\ref{4.34}), respectively. Then we can deduce the unweighted estimates from estimates (\ref{4.46}) and (\ref{4.50}) that
\begin{align}\label{4.50b}
\|v_{pxy}\|_{L^2(\Omega_\infty)}\leq C(L,\kappa)\e^{-\kappa},~~\|\p_xv_{pxy}\|_{L^2(\Omega_\infty)}\leq C(L)\e^{-1}.
\end{align}
For the corresponding weighted estimates, we recall the notations in the proof of Lemma \ref{L4.1}
\begin{align*}
&\phi(x,y):=u^0u_p-u^0_y\psi,\ \ \phi_y=u^0u_{py}-u^0_{yy}\psi\\
&\phi_x=-u^0_{xy}\psi-(v^0_p+v^1_e)u_{py}+u_{pyy}+F_p,\\
&\phi_{xx}=-u^0_{xxy}\psi+u^0_{xy}v_p-[v^0_{px}+v^1_{ex}]u_{py}+[v^0_p+v^1_e]v_{pyy}-v_{pyyy}+F_{px},
\end{align*}
where $\psi(x,y):=-\int_y^\infty u_p(x,\theta)d\theta$ and \[F_p=-u^0_{px}u^1_e-u^0_pu^1_{ex}-(v^0_p+yu^0_{px})u^0_{ez}-yu^0_{py}v^1_{ez}+E_1.\]
It is easy to see that
\begin{align*}
&\|\y^nF_p\|^2_{L^2(\Omega_\infty)}\lesssim 1+\|u_e^1\|^2_\infty+\int_0^L\sup_z|v^1_{ez}|^2\leq C,\\
&\|\y^nF_{px}\|^2_{L^2(\Omega_\infty)}\lesssim 1+\int_0^L\left(\sup_z|v^1_e|^2+\sup_z|v^1_{ez}|^2+\sup_z|v^1_{exz}|^2\right)\leq C(L,\kappa)\e^{-2\kappa},
\end{align*}
where we have used the inequality
\[|f(0)|^2=\int_0^1[(z-1)|f(z)|^2]_zdz\leq \|f\|^2_2+\|f\|_p\|f_z\|_q,~~\frac{1}{p}+\frac{1}{q}=1.\]
Then, thanks to Lemma \ref{L4.1}, there holds
\begin{align*}
\|y^n\phi_y\|_{L^2(\Omega_\infty)}&\lesssim \|\y^m\bar{u}_1\|_{H^1}+\sup_x\|\y^mv_{pyy}\|_2\leq C(L,\kappa)\e^{-\kappa},\\
\|y^n\phi_x\|_{L^2(\Omega_\infty)}&\lesssim \|\y^m\bar{u}_1\|_{H^1}+\sup_x\|\y^mv_{pyy}\|_2+\|\y^nF_p\|_{L^2(\Omega_\infty)}
\leq C(L,\kappa)\e^{-\kappa},\\
\|y^n\phi_{xx}\|_{L^2(\Omega_\infty)}&\lesssim \|\y^m\bar{u}_1\|_{H^1}+(1+\|v_e^1\|_{H^2})\sup_x\|\y^mv_{pyy}\|_2\\
&\quad+\|y^nv_{pyyy}\|_{L^2(\Omega_\infty)}+\|\y^nF_{px}\|_{L^2(\Omega_\infty)}\leq C(L,\kappa)\e^{-\kappa}.
\end{align*}
Here we again remind the readers of the fact that $y^n\leq \y^n$, for any $n\in\mathbb{N}_+.$ Further, note that
\[v_p=u^0_{px}\int_y^\infty\frac{\phi}{(u^0)^2}d\theta+u^0\int_y^\infty\left(\frac{\phi}{(u^0)^2}\right)_xd\theta.\]
Thus,
\begin{align*}
\|y^nv_{pxy}\|_{L^2(\Omega_\infty)}&\lesssim \|y^n(\phi_y,\phi_x,\phi_{xx})\|_{L^2(\Omega_\infty)}\leq C(L,\kappa)\e^{-\kappa},
\end{align*}
which, together with (\ref{4.50b}) gives $\|\y^nv_{pxy}\|_{L^2(\Omega_\infty)}\leq C(L,\kappa)\e^{-\kappa}$, and hence
\[\|\y^nv_{px}\|_{L^2(\Omega_\infty)}\leq C\|\y^{n+2}v_{pxy}\|_{L^2(\Omega_\infty)}\leq C(L,\kappa)\e^{-\kappa}.\]
In addition, applying $\p_x$ to $\phi_{xx}$ yields
\begin{align*}
\phi_{xxx}=v^0_{xxyy}\psi+u^0_{xy}v_{px}-[v^0_{pxx}+v^1_{exx}]u_{py}+[v^0_p+v^1_e]v_{pxyy}-v_{pxyyy}+F_{pxx},
\end{align*}
in which we infer that
\[\|\y^n F_{pxx}\|_{L^2(\Omega_\infty)}\lesssim 1+\e^{-\frac{1}{4}}\left(\|v^1_{ez}\|_{L^2(\Omega_0)}+\|v^1_{exz}\|_{L^2(\Omega_0)}+\|v^1_{exxz}\|_{L^2(\Omega_0)}\right)\leq C\e^{-\frac{3}{4}}.\]
Further, there holds
\begin{align*}
\|y^n\phi_{xxx}\|_{L^2(\Omega_\infty)}&\lesssim \left(1+\int_0^L\sup_z|v^1_{exx}|\right)\|\y^nv_{pyy}\|_{L^2(\Omega_\infty)}+\|\y^nv_{pxy}\|_{L^2(\Omega_\infty)}\\
&\quad+\|\y^nv_{pxyy}\|_{L^2(\Omega_\infty)}+\|y^nv_{pxyyy}\|_{L^2(\Omega_\infty)}+\|\y^nF_{pxx}\|_{L^2(\Omega_\infty)}\\
&\leq C(L)\e^{-1}.
\end{align*}
Therefore, we get
\begin{align*}
\|y^nv_{pxxy}\|_{L^2(\Omega_\infty)}\lesssim& \|y^n(\phi_y,\phi_x,\phi_{xx},\phi_{xxx})\|_{L^2(\Omega_\infty)}\leq C(L)\e^{-1}.
\end{align*}
Similarly, one can deduce the second part of (\ref{4.50a}) and hence finish the proof of it.
\endProof

Now, similar to Section 2, we cut-off the solutions from $\Omega_\infty$ to $\Omega_\e$ and prove the following
\begin{prop}\label{p3.1}
Under the assumptions in Theorem 1.1, there exists smooth functions $[u^1_p,v^1_p]$, defined in $\Omega_\e$, satisfying the following inhomogeneous system:
 \begin{align}\label{4.51}
 \begin{cases}
 u^0u^1_{px}+u^0_xu^1_p+u^0_yv^1_p+[v^0_p+v_e^1]u^1_{py}-u^1_{pyy}=R^{u,1}_p,\\
u^1_{px}+v^1_{py}=0,\\
u^1_p(0,y)=\bar{u}_1(y),u^1_p(x,0)=-u^1_e(x,0),\\
u^1_{py}(x,\frac{1}{\sqrt\e})=v^1_p(x,0)=v^1_p(x,\frac{1}{\sqrt\e})=0,
 \end{cases}
 \end{align}
 where the inhomogeneous term $R^{u,1}_p$ is a higher order term of $\sqrt\e$. In addition, there holds that
\begin{align}
&\|[u^1_p,v^1_p]\|_{L^\infty(\Omega_\e)}+\sup_{0\leq x\leq L}\|\y^nv^1_{pyy}\|_{L^2(I_\e)}+\|\y^n(v^1_{px},v^1_{pxy})\|_{L^2(\Omega_\e)}\leq C(L,\kappa)\e^{-\kappa},\label{4.52}\\
&\sup_{0\leq x\leq L}\|\y^nv^1_{pxyy}\|_{L^2(I_\e)}+\|\y^n(v^1_{pxx},v^1_{pxxy})\|_{L^2(\Omega_\e)}\leq C(L)\e^{-1},\label{4.53}
\end{align}
for any given $n\in\mathbb{N}$.
\end{prop}
\Proof
 Let $[u_p,v_p]$ be constructed in Lemma \ref{L4.3} and define that
 \begin{align*}
 &u^1_p(x,y):=\chi(\sqrt\e y)u_p(x,y)-\sqrt\e\chi^\prime(\sqrt\e y)\int_0^y u_p(x,\theta)d\theta,\\
 &v^1_p(x,y):=\chi(\sqrt\e y)v_p(x,y).
 \end{align*}
 Then it follows directly from calculation that $[u^1_p,v^1_p]$ satisfies (\ref{4.51}) with
\begin{align}\label{4.54}
R^{u,1}_p:=&\sqrt\e\chi^\prime u^0v_p+\sqrt\e\chi^\prime u^0_x\int_0^y u_pd\theta+2\sqrt\e\chi^\prime [v^0_p+v^1_e]u_p-3\sqrt\e\chi^\prime u_{py}\nonumber\\
&+\sqrt\e F_p\int_0^y\chi^\prime d\theta+\e\chi^{\prime\prime}[v^0_p+v^1_e]\int_0^y u_pd\theta-3\e\chi^{\prime\prime}u_p
-\e^{3/2}\chi^{\prime\prime\prime}\int_0^y u_pd\theta.
\end{align}
Clearly, by the estimates in Lemma \ref{L4.3}, we get
\[\Big|\sqrt\e\chi^\prime(\sqrt\e y)\int_0^y u_p(x,\theta)d\theta\Big|\leq \sqrt\e y|\chi^\prime(\sqrt\e y)|\|u_p\|_\infty\leq C(L,\kappa)\e^{-\kappa}.\]
Hence, (\ref{4.52}) and (\ref{4.53}) follows from Lemma \ref{L4.3}. The proof of this lemma is completed.
\endProof

\begin{cor}\label{c4.6}
Assume that $[u^1_p,v^1_p,p^1_p]$ is the solution to (\ref{4.51}), where $p^1_p$ is an absolute constant. Then, for any $\kappa>0$ small sufficiently, there holds that
\begin{align}\label{4.55}
\|\tilde{R}^u_1\|_{L^2(\Omega_\e)}+\|R^{u,1}_p\|_{L^2(\Omega_\e)}\leq C(L,\kappa)\e^{\frac{1}{4}-\kappa},\ \ \|p^2_{px}\|_{L^2(\Omega_\e)}\leq C\e^{-\frac{1}{4}}
\end{align}
where $\tilde{R}^u_1,p^2_p$ are defined as in (\ref{4.0}) and (\ref{3.7}),respectively. is a small enough constant.
\end{cor}
\Proof
Thanks to Lemma \ref{L4.3}, there holds that
\begin{align}
\|\tilde{R}^u_1\|^2_{L^2(\Omega_\e)}\lesssim&\e^{1/2}\|\y^nv^1_{pyy}\|^2_{L^2(\Omega_\e)}\|u^0_{ez}\|^2_{L^2(0,1)}
+\e^{1/2}\|v^1_p\|^2_\infty\|u^0_{ez}\|^2_{L^2(0,1)}\nonumber\\
\leq& C(L,\kappa)\e^{\frac{1}{2}-2\kappa},\label{4.56}\\
\|R^{u,1}_p\|^2_{L^2(\Omega_\e)}\lesssim& \sqrt\e L\|v_p\|^2_\infty
+\e\|u_p\|^2_\infty\iint_{\Omega_\e}\y^2|u^0_{px}|^2\nonumber\\
&+\e\|u_p\|^2_\infty\iint_{\Omega_\e}|v^0_p+v^1_e|^2+\e\iint_{\Omega_\e}|u_{py}|^2+\e\iint_{\Omega_\e}\y^2|F_p|^2\nonumber\\
&+\e\|u_p\|^2_\infty\iint_{\Omega_\e}|v^0_p+v^1_e|^2+\e^{3/2}\|u_p\|^2_\infty\nonumber\\
\leq &C(L,\kappa)\e^{\frac{1}{2}-2\kappa}.\label{4.561}
\end{align}
In view of equation (\ref{3.7}), we get
\begin{align}\label{4.57}
p^2_{px}=\int_y^{1/\sqrt\e}\left[(u^0_e+u^0_p)v^0_{pxx}+u^0_pv^1_{exx}+(v^0_p+v^1_e)v^0_{pxy}-v^0_{pxyy}\right](x,\theta)d\theta.
\end{align}
Note that
\begin{align*}
\int_y^{1/\sqrt\e}[u^0_e+u^0_p]v^0_{pxx}&\leq C\y^{-n+1}\|u^0_e+u^0_p\|_\infty\|\y^nv^0_{pxx}\|_2,\\
\int_y^{1/\sqrt\e}u^0_pv^1_{exx}&\leq C\y^{-n+2}\|\y^nu^0_{py}\|_2\|v^1_{exx}(x,\sqrt\e\cdot)\|_2,\\
\int_y^{1/\sqrt\e}(v^0_p+v^1_e)v^0_{pxy}&\leq C\y^{-n+1}\|u^0_e+u^0_p\|_\infty\|\y^nv^0_{pxy}\|_2,\\
\int_y^{1/\sqrt\e}v^0_{pxyy}&\leq C\y^{-n+1}\|\y^nv^0_{pxyy}\|_2.
\end{align*}
Hence, taking $n\geq 3$, we have
\begin{align}\label{4.58}
\|p^2_{px}\|_{L^2(\Omega_\e)}\leq C\e^{-1/4}.
\end{align}
The proof of this corollary is completed.
\endProof

With the zeroth order Prandtl profiles, first order Euler correctors and the first order Prandtl corrector in hand, together with those various estimates on the approximate solution, we are able to give the error estimates as follows.
\begin{prop}\label{p4.1}
Under the same assumptions as in Theorem 1.1, there exists approximate solutions $[u_{app},v_{app},p_{app}]$ such that
\begin{align}\label{400}
\|R^u_{app}\|_{L^2(\Omega_\e)}+\sqrt\e\|R^v_{app}\|_{L^2(\Omega_\e)}\leq C(L,\kappa)\e^{\frac{3}{4}-\kappa},
\end{align}
where $C$ depends on initial data and $L,\kappa$.
\end{prop}
\Proof
Collecting errors from $R^u_0$ in (\ref{2.3}), $R^{u,0}_p$ in (\ref{2.012}), $R^u_1$ in (\ref{3.23}), $\tilde{R}^u_1$ in (\ref{4.0}), $R^{u,1}_p$ in (\ref{4.54}), and the remaining $\e$-order terms in $R^u_{app}$, we yield
\begin{align*}
R^u_{app}=&E_0-\e u^0_{ezz}+\e E_2+\sqrt\e R^u_1+\sqrt\e\tilde{R}^u_1+\sqrt\e R^{u,1}_p\\
&+\e\left[(u^1_e+u^1_p)\p_x+v^1_p\p_y\right](u^1_e+u^1_p)+\e p^2_{px}-\e\p_x^2\left[u^0_p+\sqrt\e(u^1_e+u^1_p)\right]
\end{align*}
In view of (\ref{3.25}), (\ref{3.23}), (\ref{3.24}), and (\ref{4.55}), we immediately get
\begin{align}\label{401}
\|E_0-\e u^0_{ezz}+\sqrt\e R^{u}_1+\e\tilde{R}^u_1+\sqrt\e R^{u,1}_p+\e p^2_{px}\|_{L^2(\Omega_\e)}\leq C(L,\kappa)\e^{\frac{3}{4}-\kappa}.
\end{align}

Similarly, using the estimates for $[u^1_e,v^1_e]$, $[u^1_p,v^1_p]$, and $u^0_p$, we have
\begin{align*}
\e\|(u^1_e+u^1_p)\p_x(u^1_e+u^1_p)\|_{L^2(\Omega_\e)}&\leq \e\left(\|u^1_e\|_{\infty}+\|u^1_p\|_{\infty}\right)\left(\|u^1_{ex}\|_2+\|v^1_{py}\|_2\right)\nonumber\\
&\leq C(L,\kappa)\e^{1-\kappa},\\
\e\|v^1_p\p_y(u^1_e+u^1_p)\|_{L^2(\Omega_\e)}&\leq \e\|v^1_p\|_{\infty}\left(\sqrt\e\|u^1_{ez}\|_{2}+\|u^1_{py}\|_{2}\right)\leq C(L,\kappa)\e^{1-\kappa},\\
\e\|\p_x^2(u^0_p+\sqrt\e u^1_e+\sqrt\e u^1_p)\|_{L^2(\Omega_\e)}&\leq \e\|u^0_{pxx}\|_{2}+\e^{\frac{3}{2}}\left(\|v^1_{exz}\|_{2}+\|v^1_{pxy}\|_{2}\right)\leq C\e^{\frac{3}{2}-\kappa},
\end{align*}
where we have used Proposition \ref{p3.1} with $\kappa/2$.

In addition, recalling the definition of $E_2$ in (\ref{2.012}), one has
\begin{align*}
\|E_2\|_{L^2(\Omega_\e)}\lesssim& \|\y u^\infty_p\|_{2}\|v_p^\infty\|_{\infty}+\|u^\infty_p\|_{2}+\|\y u^\infty_p\|_{2}\|v^\infty_p\|_2\\
&+\|\chi^{\prime\prime}\|_2\left(\|v^\infty_p\|_{2}+\|v^\infty_{py}\|_{2}\right)+\|\y u^\infty_p\|_{2}\|\chi^{\prime\prime\prime}\|_2\\
\leq& C(L).
\end{align*}
Thus, we deduce that
\[\|R^{u}_{app}\|_{L^2(\Omega_\e)}\leq C(L,\kappa)\e^{\frac{3}{4}-\kappa}.\]

Next, we will give estimates for $R^v_{app}$. Recalling the remaining terms of $R^v_0$ in (\ref{3.701}) and the definition of $R^v_{app}$, we infer that
\begin{align*}
R^v_{app}=&R^v_0+\sqrt\e\left[(u^0_e+u^0_p+\sqrt\e[u^1_e+u^1_p])\p_x+([v^0_p+v^1_e]+\sqrt\e v^1_p)\p_y\right]v^1_p\\
&+\sqrt\e\left[(u^1_e+u^1_p)\p_x+v^1_p\p_y\right](v^0_p+v^1_e)-\sqrt\e v^1_{pyy}-\e\p_x^2\left[(v^0_p+v^1_e)+\sqrt\e v^1_p\right].
\end{align*}
Similarly as above, we have that
\begin{align*}
&\sqrt\e\left\|\left[(u^0_e+u^0_p+\sqrt\e[u^1_e+u^1_p])\p_x+([v^0_p+v^1_e]+\sqrt\e v^1_p)\p_y\right]v^1_p\right\|_{L^2(\Omega_\e)}\\
\leq& C\sqrt\e\|[u^0_e,u^0_p,v^0_p,u^1_e,v^1_e]\|_{\infty}\left(\|v^1_{px}\|_{L^2}+\|v^1_{py}\|_{2}\right)
\leq C(L,\kappa)\e^{\frac{1}{4}-\kappa},
\end{align*}
and that
\begin{align*}
&\sqrt\e\left\|\left[(u^1_e+u^1_p)\p_x+v^1_p\p_y\right](v^0_p+v^1_e)\right\|_{L^2(\Omega_\e)}\\
\leq& C\sqrt\e\|[u^1_e,u^1_p,v^1_p]\|_{\infty}\left(\|v^0_{px}+v^1_{ex}\|_{2}+\|v^0_{py}+\sqrt\e v^1_{ez}\|_{2}\right)\\
\leq& C(L,\kappa)\e^{\frac{1}{4}-\kappa}.
\end{align*}
Moreover, one has
\[\|\sqrt\e v^1_{pyy}-\e\p_x^2(v^0_p+v^1_e+\sqrt\e v^1_p)\|_{L^2(\Omega_\e)}\leq C(L)\e^{\frac{1}{4}-\kappa}.\]
Putting these estimates together and using the estimate (\ref{3.24}) yield
\[\|R^v_{app}\|_{L^2(\Omega_\e)}\leq C(L,\kappa)\e^{\frac{1}{4}-\kappa},\]
which completes the proof of this Proposition.
\endProof

\section{The existence of remainder solutions}
Now we are on the final step to prove the main theorem: the existence of the remainder solutions. Denote that
\begin{align}\label{5.1}
u_s(x,y):=u^0_e(\sqrt\e y)+u^0_p(x,y)+\sqrt\e u^1_e(x,\sqrt\e),\ \ v_s(x,y):=v^0_p(x,y)+v^1_e(x,\sqrt\e y).
\end{align}
Then the remainder solutions $[u^\e,v^\e,p^\e]$ solves
\begin{equation}\label{5.2}
\begin{cases}
u_su^\e_x+u^\e u_{sx}+v_su^\e_y+v^\e u_{sy}+p^\e_x-\Delta_\e u^\e=R_1(u^\e,v^\e),\\
u_sv^\e_x+u^\e v_{sx}+v_sv^\e_y+v^\e v_{sy}+p^\e_y/\e-\Delta_\e v^\e=R_2(u^\e,v^\e),\\
u^\e_x+v^\e_y=0,
\end{cases}
\end{equation}
where
\begin{align*}
&R_1(u^\e,v^\e):=\e^{-\gamma-\frac{1}{2}}R^u_{app}-\sqrt\e\left[(u^1_p+\e^\gamma u^\e)u^\e_x+u^\e u^1_{px}+(v^1_p+\e^\gamma v^\e)u^\e_y+v^\e u^1_{py}\right],\\
&R_2(u^\e,v^\e):=\e^{-\gamma-\frac{1}{2}}R^v_{app}-\sqrt\e\left[(u^1_p+\e^\gamma u^\e)v^\e_x+u^\e v^1_{px}+(v^1_p+\e^\gamma v^\e)v^\e_y+v^\e v^1_{py}\right].
\end{align*}
The errors $R^u_{app}$ and $R^v_{app}$ in $R_1$ and $R_2$ have been estimated in Proposition \ref{p4.1}. It should be noted that, since $\min_y\{u^0_e(\sqrt\e y)+\bar{u}_0(y)\}>0$ and $\|u^1_e\|_\infty\leq C$, the known function $u_s$ in (\ref{5.2}) is strictly positive as $\e$ and $L$ small sufficiently. This is very important in using the positivity as is done in (\ref{3.14}) and (\ref{3.15}).

Before begining to prove the existence of the remainder $[u^\e,v^\e,p^\e]$, we first give the following two Propositions, the proof of which are stated in Section 3 and Section 4 of \cite{GN14}, respectively, and hence we omit the detail here.

The first proposition gives the linear stability estimates for (\ref{5.2}):
\begin{prop}\label{p5.1}
For any given $f,g\in L^2(\Omega_\e)$, there exists some positive number $L$ such that the linear problem
\begin{equation}\label{5.3}
\begin{cases}
u_su_x+uu_{sx}+v_su_y+vu_{sy}+p_x-\Delta_\e u=f,\\
u_sv_x+uv_{sx}+v_sv_y+vv_{sy}+p_y/\e-\Delta_\e v=g,\\
u_x+v_y=0,
\end{cases}
\ \textrm{in}\ \Omega_\e,
\end{equation}
together with boundary conditions
\begin{align}\label{5.4}
\begin{cases}
[u,v]_{y=0}=0,\ \ [u_y,v]_{y=\frac{1}{\sqrt\e}}=0,\\
[u,v]_{x=0}=0,\ \ [p-2\e u_x,u_y+\e v_x]_{x=L}=0,
\end{cases}
\end{align}
has an unique solution $[u,v,p]$ defined on $\Omega_\e$. In addition, there holds
\begin{align}\label{5.5}
\|\nabla_\e u\|_{L^2(\Omega_\e)}+\|\nabla_\e v\|_{L^2(\Omega_\e)}\lesssim\|f\|_{L^2(\Omega_\e)}+\sqrt\e\|g\|_{L^2(\Omega_\e)}.
\end{align}
\end{prop}

The second one provides $L^\infty$ estimates of the solution to the corresponding Stokes problem:
\begin{prop}\label{p5.2}
For any given $f,g\in L^2(\Omega_\e)$, consider the incompressible Stokes equation
\begin{equation}\label{5.6}
\begin{cases}
-\Delta_\e u+p_x=f,\\
-\Delta_\e v+p_y/\e=g,\\
u_x+v_y=0,
\end{cases}
\ \textrm{in}\ \Omega_\e,
\end{equation}
together with the same boundary conditions as in (\ref{5.4}). Then, for any $\gamma>0$, there holds that
\begin{align}\label{5.7}
&\|u\|_{L^\infty(\Omega_\e)}+\sqrt\e\|v\|_{L^\infty(\Omega_\e)}\nonumber\\
\lesssim& C_{\gamma,L}\e^{-\frac{\gamma}{4}}\left(\|\nabla_\e u\|_{L^2(\Omega_\e)}+\|\nabla_\e v\|_{L^2(\Omega_\e)}+\|f\|_{L^2(\Omega_\e)}+\sqrt\e\|g\|_{L^2(\Omega_\e)}\right),
\end{align}
for some constant $C_{\gamma,L}$ depending only on $\gamma$ and $L$.
\end{prop}

\ProofTheorem ~~~With these two propositions in hand, we are able to apply the standard contraction mapping principle for the existence of solutions to the nonlinear problem, which is consisted of several steps.

Step 1. We introduce the function space $\mathscr{X}$ endowed with the norm:
\begin{align}\label{5.8}
\|[u^\e,v^\e]\|_{\mathscr{X}}:=\|\nabla_\e u^\e\|_{L^2(\Omega_\e)}+\|\nabla_\e v^\e\|_{L^2(\Omega_\e)}+\|u^\e\|_{L^\infty(\Omega_\e)}+\sqrt\e\|v^\e\|_{L^\infty(\Omega_\e)},
\end{align}
where $\nabla_\e:=\p_y+\sqrt\e\p_x.$ And, we choose the following subspace of $\mathscr{X}$ with $K$ to be determined:
\[\mathscr{X}_K:=\left\{[u^\e,v^\e]\in\mathscr{X}\big|\|[u^\e,v^\e]\|_{\mathscr{X}}\leq K\right\}.\]

Step 2. For each $[\bar{u}^\e,\bar{v}^\e]\in\mathscr{X}_K$, we solve the corresponding linearized problem for $[u^\e,v^\e]$:
\begin{equation}\label{5.9}
\begin{cases}
u_su^\e_x+u^\e u_{sx}+v_su^\e_y+v^\e u_{sy}+p^\e_x-\Delta_\e u^\e=R_1(\bar{u}^\e,\bar{v}^\e),\\
u_sv^\e_x+u^\e v_{sx}+v_sv^\e_y+v^\e v_{sy}+p^\e_y/\e-\Delta_\e v^\e=R_2(\bar{u}^\e,\bar{v}^\e),\\
u^\e_x+v^\e_y=0,
\end{cases}
\end{equation}
equips with the same boundary conditions as (\ref{5.4}). Then, by Proposition \ref{p5.1}, there exists an unique strong solution $[u^\e,v^\e,p^\e]$ satisfying that
\begin{align}\label{5.10}
\|\nabla_\e u^\e\|_{L^2(\Omega_\e)}+\|\nabla_\e v^\e\|_{L^2(\Omega_\e)}\leq \|R_1(\bar{u}^\e,\bar{v}^\e)\|_{L^2(\Omega_\e)}+\sqrt\e\|R_2(\bar{u}^\e,\bar{v}^\e)\|_{L^2(\Omega_\e)}.
\end{align}
Now, we give estimates for $R_1,R_2.$ In view of Proposition \ref{p4.1}, for any $\kappa>0$, it follows that
\begin{align}\label{5.11}
\e^{-\gamma-\frac{1}{2}}\left[\|R^u_{app}\|_{L^2(\Omega_\e)}+\sqrt\e\|R^v_{app}\|_{L^2(\Omega_\e)}\right]\leq C(L,\kappa)\e^{\frac{1}{4}-\gamma-\kappa}.
\end{align}
In addition, using the estimates for $[u^1_p,v^1_p]$ in Proposition \ref{p3.1} and the divergence-free condition $u_x^\e+v^\e_y=0$, we infer that
\begin{align}\label{5.12}
&\sqrt\e\|(u^1_p+\e^\gamma\bar{u}^\e)\bar{u}^\e_x+(v^1_p+\e^\gamma\bar{v}^\e)\bar{u}^\e_y\|_{L^2(\Omega_\e)}\nonumber\\
\leq& \sqrt\e\left[(\|u^1_p\|_{\infty}+\e^\gamma\|\bar{u}^\e\|_{\infty})\|\bar{v}^\e_y\|_{2}
+(\|v^1_p\|_{\infty}+\e^\gamma\|\bar{v}^\e\|_{\infty})\|\bar{u}^\e_y\|_{2}\right]\nonumber\\
\leq &\sqrt\e\|[u^1_p,v^1_p]\|_{\infty}\|[\bar{u}^\e,\bar{v}^\e]\|_{\mathscr{X}}
+\e^\gamma\|[\bar{u}^\e,\bar{v}^\e]\|^2_{\mathscr{X}}\nonumber\\
\leq& C(L,\kappa)\e^{\frac{1}{2}-\kappa}K+\e^\gamma K^2,
\end{align}
and that
\begin{align}\label{5.13}
\sqrt\e\|\bar{u}^\e u^1_{px}+\bar{v}^\e u^1_{py}\|_{L^2(\Omega_\e)}&\leq C\sqrt\e\left[\|\bar{u}_y^\e\|_{2}\sup_x\|\y^nu^1_{px}\|_2
+\|\bar{v}_y^\e\|_{2}\sup_x\|\y^nu^1_{py}\|_2\right]\nonumber\\
&\leq C(L,\kappa)\e^{\frac{1}{2}-\kappa}\|[\bar{u}^\e,\bar{v}^\e]\|_{\mathscr{X}}\leq C(L,\kappa)\e^{\frac{1}{2}-\kappa}K,
\end{align}
in which we have used the fact that $|[u^\e,v^\e]|\leq \sqrt y\|[u^\e_y,v^\e_y]\|_2$.

Similarly, for the term in $R_2$, there holds that
\begin{align}\label{5.14}
&\sqrt\e\|(u^1_p+\e^\gamma\bar{u}^\e)\bar{v}^\e_x+(v^1_p+\e^\gamma\bar{v}^\e)\bar{v}^\e_y\|_{L^2(\Omega_\e)}\nonumber\\
\leq& (\|u^1_p\|_{\infty}+\e^\gamma\|\bar{u}^\e\|_{\infty})\|\sqrt\e\bar{v}^\e_x\|_{2}
+\sqrt\e(\|v^1_p\|_{\infty}+\e^\gamma\|\bar{v}^\e\|_{\infty})\|\bar{v}^\e_y\|_{2}\nonumber\\
\leq&C\|[u^1_p,v^1_p]\|_{\infty}\|[\bar{u}^\e,\bar{v}^\e]\|_{\mathscr{X}}
+\e^\gamma\|[\bar{u}^\e,\bar{v}^\e]\|^2_{\mathscr{X}}\nonumber\\
\leq& C(L,\kappa)\e^{-\kappa}K+\e^\gamma K^2,
\end{align}
and that
\begin{align}\label{5.15}
&\sqrt\e\|\bar{u}^\e v^1_{px}+\bar{v}^\e v^1_{py}\|_{L^2(\Omega_\e)}\nonumber\\
\leq&C\sqrt\e\left[\|\bar{u}^\e\|_{L^\infty}\|v^1_{px}\|_{2}
+\|\bar{v}_y^\e\|_{2}\sup_x\|\y^nv^1_{pyy}\|_2\right]\nonumber\\
\leq& C(L,\kappa)\e^{\frac{1}{2}-\kappa}\|[\bar{u}^\e,\bar{v}^\e]\|_{\mathscr{X}}\leq C(L,\kappa)\e^{\frac{1}{2}-\kappa}K.
\end{align}
where the estimate $\|v^1_{px}\|_{L^2(\Omega_\e)}\leq C\e^{-\kappa}$ has been used.

In conclusion, we yield
\begin{align}\label{5.16}
&\|R_1(\bar{u}^\e,\bar{v}^\e)\|_{L^2(\Omega_\e)}+\sqrt\e\|R_2(\bar{u}^\e,\bar{v}^\e)\|_{L^2(\Omega_\e)}\nonumber\\
&\leq C(L,\kappa)\e^{\frac{1}{4}-\kappa-\gamma}+C(L,\kappa)\e^{\frac{1}{2}-\kappa}K+\e^\gamma K^2,
\end{align}
which implies the estimate for the gradient of $[u^\e,v^\e]$:
\begin{align}\label{5.17}
&\|\nabla_\e u^\e\|_{L^2(\Omega_\e)}+\|\nabla_\e v^\e\|_{L^2(\Omega_\e)}\nonumber\\
\leq &C(u_s,v_s,L,\kappa)\e^{\frac{1}{4}-\kappa-\gamma}+C(u_s,v_s,L,\kappa)\e^{\frac{1}{2}-\kappa}K+C(u_s,v_s)\e^\gamma K^2.
\end{align}

It remains to estimate the $L^\infty$ norm for $[u^\e,v^\e]$. Recalling Proposition \ref{p5.2} with
\begin{align*}
&f:=R_1-u_su_x^\e-u_{sx}u^\e-v_s\bar{u}^\e_y-v^\e u_{sy},\\
&g:=R_2-u_sv_x^\e-v_{sx}u^\e-v_s\bar{v}^\e_y-v^\e v_{sy},
\end{align*}
it follows from (\ref{5.7}) that
\begin{align}\label{5.18}
&\|u^\e\|_{L^\infty(\Omega_\e)}+\sqrt\e\|v^\e\|_{L^\infty(\Omega_\e)}\nonumber\\
\leq& C_{\gamma,L}\e^{-\frac{\gamma}{4}}\left(\|\nabla_\e u^\e\|_2+\|\nabla_\e v^\e\|_2
+\|R_1-u_su_x^\e-u_{sx}u^\e-v_s\bar{u}^\e_y-v^\e u_{sy}\|_2\right)\nonumber\\
&+C_{\gamma,L}\e^{\frac{1}{2}-\frac{\gamma}{4}}\|R_2-u_sv_x^\e-v_{sx}u^\e-v_s\bar{v}^\e_y-v^\e v_{sy}\|_2.
\end{align}
Since that (\ref{5.17}) and (\ref{5.16}) have give the desired estimate for $[\nabla_\e u^\e,\nabla_\e v^\e]$ and $[R_1,R_2]$, respectively, it remains to estimate the rest terms with respect to $[u_s,v_s]$.  Indeed, note that
\begin{align*}
\sup_x\|\sqrt yu_{sx}\|_{L^2(I_\e)}&\leq \sup_x\|\y u^0_{px}\|_2+\sup_x\|v^1_{ez}\|_2\\
&\leq \sup_x\|\y u^0_{px}\|_2+\|v^1_{ezx}\|_2+\|\p_zV_{b0}\|_2,\\
\sup_x\|\sqrt yu_{sy}\|_{L^2(I_\e)}&\leq \|u^0_{ez}\|_2+\sup_x\|\y u^0_{py}\|_2+\sup_x\|u^1_{ez}\|_2\\
&\leq \|u^0_{ez}\|_2+\sup_x\|\y u^0_{py}\|_2+\|u^1_{bz}\|_2+\|v^1_{ezz}\|_2,\\
\|[u_s,v_s]\|_{L^\infty(\Omega_\e)}&\leq \|u^0_e\|_\infty+\|[u^0_p,v^0_p]\|_\infty+\|[u^1_e,v^1_e]\|_\infty,
\end{align*}
we have
\begin{align*}
\|u_{sx}u^\e+u_{sy}v^\e\|_{L^2(\Omega_\e)}&\leq \left(\|u^\e_y\|_2\sup_x\|\sqrt yu_{sx}\|_2
+\|v^\e_y\|_2\sup_x\|\sqrt yu_{sy}\|_2\right)\\
&\leq C\|[\nabla_\e u^\e,\nabla_\e v^\e]\|_2,\\
\|u_su^\e_x+v_su^\e_y\|_{L^2(\Omega_\e)}&\leq \|[u_s,v_s]\|_{\infty}\|[u_y^\e,v^\e_y]\|_{2}\leq C\|[\nabla_\e u^\e,\nabla_\e v^\e]\|_{2}.
\end{align*}
Similarly, there holds that
\begin{align*}
\|v_{sx}u^\e+v_{sy}v^\e\|_{L^2(\Omega_\e)}&\leq \left(\|u^\e_y\|_{2}\sup_x\|\sqrt yv_{sx}\|_2
+\|v^\e_y\|_{2}\sup_x\|\sqrt yu_{sx}\|_2\right)\\
&\leq C\e^{-\frac{1}{2}}\|[\nabla_\e u^\e,\nabla_\e v^\e]\|_{2},\\
\|u_sv^\e_x+v_sv^\e_y\|_{L^2(\Omega_\e)}&\leq \|[u_s,v_s]\|_{\infty}\|[v_x^\e,v^\e_y]\|_{2}\leq C\e^{-\frac{1}{2}}\|[\nabla_\e u^\e,\nabla_\e v^\e]\|_{2},
\end{align*}
in which the following estimate has been used:
\begin{align*}
\sup_x\|\sqrt yv_{sx}\|_{L^2(I_\e)}&\leq \sup_x\|\y v^0_{px}\|_2+\e^{-\frac{1}{2}}\sup_x\|v^1_{ex}(x,\cdot)\|_2\\
&\leq \sup_x\|\y v^0_{px}\|_2+\e^{-\frac{1}{2}}\|v^1_{ex}\|_2\|v^1_{exx}\|_2.
\end{align*}
Substituting these estimates together with (\ref{5.17}), (\ref{5.16}) into (\ref{5.18}) then gives
\begin{align}\label{5.19}
&\|u^\e\|_{L^\infty(\Omega_\e)}+\sqrt\e\|v^\e\|_{L^\infty(\Omega_\e)}\nonumber\\
\leq& C(u_s,v_s,L,\kappa)\e^{\frac{1}{4}-\kappa-\frac{5\gamma}{4}}+C(u_s,v_s,L,\kappa)\e^{\frac{1}{2}-\kappa-\frac{\gamma}{4}}K
+C(u_s,v_s)\e^\frac{3\gamma}{4}K^2.
\end{align}

Now, adding up (\ref{5.18}) and (\ref{5.19}), noting that $\frac{1}{4}-\kappa-\frac{5\gamma}{4}\geq 0$ and $\e\ll 1$, we get
\begin{align}
\|[u^\e,v^\e]\|_{\mathscr{X}}\leq C(u_s,v_s,L,\kappa)+C(u_s,v_s,L,\kappa)\e^\frac{1}{4}K
+C(u_s,v_s)\e^\frac{3\gamma}{4}K^2.
\end{align}
Then, we take $K:=C(u_s,v_s,L,\kappa)+1$ and hence $\|[u^\e,v^\e]\|_{\mathscr{X}}\leq K$, for any small $\e$ so that
\[C(u_s,v_s,L,\kappa)\e^\frac{1}{4}K+C(u_s,v_s)\e^\frac{3\gamma}{4}K^2\leq 1.\]
This proves that the operator $\mathcal{M}:[\bar{u}^\e,\bar{v}^\e]\mapsto[u^\e,v^\e]$ maps $\mathscr{X}_K$ into itself.

Step 3. In order to apply the contraction mapping theorem, it remains to prove that the operator $\mathcal{M}$ is a contractive mapping. Indeed, for any two pairs $[\bar{u}^\e_1,\bar{v}^\e_1]$ and $[\bar{u}^\e_2,\bar{v}^\e_2]$ in $\mathscr{X}_K$, it follows from the similar approach that
\[\|[u^\e_1-u^\e_2,v^\e_1-v^\e_2]\|_{\mathscr{X}}
\leq C(u_s,v_s,L,\kappa)(\e^{\frac{1}{2}-\kappa-\frac{\gamma}{4}}+\e^\frac{3\gamma}{4}K)
\|[\bar{u}^\e_1-\bar{u}^\e_2,\bar{v}^\e_1-\bar{v}^\e_2\|_{\mathscr{X}},\]
which at once implies the contraction of $\mathcal{M}$, for any $\e$ small sufficiently.

This proves the existence of the unique solution to (\ref{5.2}) via standard contraction mapping theorem and hence completes the proof of the Theorem \ref{thm1}.
\endProofTheorem
\acknowledgment
Ding's research is supported by the National Natural Science Foundation of China (No.11371152, No.11571117, No.11771155 and No.11871005) and Guangdong Provincial Natural Science Foundation (No.2017A030313003).

\end{document}